\documentstyle[12pt]{article} \openup 7pt 
5pt \def\a{\alpha} \def\b{\beta} \def\e{\epsilon} \def\vs{\vspace*} \def
\O{\Omega}  \def \Z{\hbox{$Z$\hskip
-5.2pt $Z$}}    \def \C{\hbox{$C$\hskip -5pt \vrule height 6pt
depth 0pt \hskip 6pt}} \def\qed{\hfill \hfill \ifhmode\unskip\nobreak\fi
\ifmmode\ifinner\else{}\hskip2pt \fi\fi \hbox{\hskip1pt\vrule width4pt
height6pt depth1.5pt\hskip 1 pt}}  \def\lk#1#2{\mbox{$#1{}_{^{
\dis\sim}_{\ssc\Si}}#2$}} \def\d{\delta} \def\D{\Delta} \def\g{\gamma}
\def\G{\Gamma} \def\l{\lambda} \def\L{\Lambda} \def\o{\omega}  \def\Si{\Sigma} \def\si{\sigma} \def\v{\varphi} \def\u{\upsilon}
\def\VBL{\overline{V}\!({\sc\!}\L{\sc\!})} \def\sc{\scriptstyle} \def
\ssc{\scriptscriptstyle} \def\dis{\displaystyle} \def\cl{\centerline}
  \def\nl{\newline} \def\ol
{\overline}  \def\wt{\widetilde} \def\wh{\widehat}
\def\rar{\rightarrow}    \def\lra{\leftrightarrow} \def\Lra{\Leftrightarrow}
\def\bs{\backslash} \def\hs{\hspace*} \def\rb{\raisebox} \def\ra{\rangle
} \def\la{\langle} \def\ni{\noindent} \def\hi{\hangindent} \def\ha
{\hangafter} \def\pr#1{\mbox{\rb{-4pt}{$^{\dis\,\ \Pi\ }_{\sc #1}$}}}
 \def\zg{\ge} \def\zi{\in} \def\zl{\le}  \def\zn{\ne
}  \def\SUM{{\sc\sum}} \def\es{\par\vs{-4pt}\ \par\ni} 
\begin
{document}  \def\ABS{To study finite-dimensional modules
of the Lie superalgebras, Kac introduced the Kac-modules $\VBL$ and
divided them into typical or atypical modules according as they are
simple or not. For $\L$ being atypical, Hughes {\sl et al} have an
algorithm to determine all the composition factors of a Kac-module;
they conjectured that there exists a bijection between the composition
factors of a Kac-module and the so-called permissible codes. The aim of
this paper is to prove this conjecture. We gives a partial proof here,
i.e., to any unlinked code, by constructing explicitly the primitive
vector, we prove that there corresponds a composition factor of the
Kac-module. We will give a full proof of the conjecture in another
paper.} \def\KW{Kac-module, atypical, composition factor, primitive
weight, code.}  \par\ \par \par
\ \par \cl{\bf PRIMITIVE VECTORS OF KAC-MODULES OF} \cl{\bf THE LIE
SUPERALGEBRAS $sl(m/n)$} \par\ \par\cl{Yucai Su$^{\,a,}$\footnote{$^)$
Partly supported by a grant from Shanghai Jiaotong University}$^)$, \
J. W. B. Hughes$^{\,b}$ \ and \ R. C. King$^{\,c}$}\par\ni\hs{3ex} $^{a)
}$ Department of Applied Mathematics, Shanghai Jiaotong University,
China\nl\hs{3ex} $^{b)}$ School of Mathematical Sciences, Queen Mary
and Westfield College, London, U.K.\nl\hs{3ex} $^{c)}$ Faculty of
Mathematical Studies, University of Southampton, U.K.\par\ \par \ni{\bf
ABSTRACT.}\ \ABS\par\ni{\bf KEYWORDS:}\ \KW \vs{-5pt}\par\ \es{\bf I.
INTRODUCTION}\par\ni In the classification of finite-dimensional
modules of the basic classical Lie superalgebras,$^{3-6,10}$ Kac
distinguished between typical and atypical modules. He also introduced
now the so-called Kac-module $\VBL$, which was shown to be simple if
and only if $\L$ is typical. For $\L$ being atypical, the problem of
the structure of $\VBL$, or equivalently, the character of the simple
module $V(\L)$, has been the subject of intensive study.$^{2,11,14-16}$
More generally, the problem of classifying indecomposable modules has
received much attention in the literature.$^{1,8,9,11,12}$ Kac obtained
a character formula for typical modules.$^5$ The problem for atypical $
sl(m/n)$-modules has seemed to be difficult, though several partial
solutions have been achieved.$^{2,15,16}$\par Only recently Serganova$^
{\,11}$ found a solution for the characters of simple $gl(m/n)
$-modules, who \,described \,the \,multiplicities \,$a_{\ssc\L\Si}$
\,of \,composition \,factors \,$V(\Si)$ \,of \,$\VBL$ \,in \,terms \,of
Kazhdan-Lusztig polynomials. However, Serganova's algorithm of
describing $a_{\ssc\L\Si}$ turns out to be rather complicated. The
structure of $\VBL$ is still not so apparent to readers.\vs{-2pt}
Hughes {\sl et al}$^{\,2}$ derived an\vs{-2pt} algorithm to determine
all the composition factors of $sl(m/n)$-Kac-modules $\VBL$. They
conjectured\vs{-2pt} that there exists a bijection between the
composition factors of $\VBL$ and the permissible codes (Definition
3.9). This conjecture\vs{-2pt} clearly describes the structure of $\VBL
$. \vs{-2pt}The aim of the present paper and the forthcoming paper$^{
\,13}$ is to prove this conjecture. In this paper, we prove that to any
unlinked code, there corresponds a composition factor of $\VBL$, by
constructing explicitly a primitive vector corresponding to the
unlinked code (Theorems 6.6\&6.12). Then in {\sl Ref.} 13, we will give
a full proof of the conjecture and point out that the proof of the
conjecture will result in the proofs of some other conjectures. \es{\bf
II. THE LIE SUPERALGEBRA $sl(m+1/n+1)$}\vs{-2pt}\par\ni \vs{-1pt}Denote
$G$=$sl(m$+$1/n$+$1)$ the set of ($m$+$n$+2)$\times$($m$+$n$+2)
matrices $x$=$(^A_C\,^B_D)$ of zero supertrace $str(x)$=$tr(A)${}$-${}$
tr(D)$=0, \vs{-1pt}where $A,B,C,D$ are ($m$+1)$\times$($m$+1),($m$+1)$
\times$($n$+1),($n$+1)$\times$($m$+1), ($n$+1)$\times$($n$+1) matrices
respectively. Let $G_{\ol0}$=$\{(^A_0\,^0_D)\},$ $G_{\ol1}$=$\{(^0_C\,^
B_0)\}$, then $G$=$G_{\ol0}${}$\oplus${}$G_{\ol1}$ as a\vs{-1pt} $\Z_2$
=$\Z$/2$\Z$ graded space, is a Lie superalgebra with\vs{-1pt} bracket:
[$x$,$y$]=$xy${}$-$($-$1)$^{\xi\eta}yx$ for $x${}$\in${}$G_\xi$, $y${}$
\in${}$G_\eta,\xi$,$\eta${}$\in${}$\Z_2$ such that $G_{\ol0}${}$\cong
${}$sl(m$+1)$\oplus${}$\C${}$\oplus${}$sl(n$+1) is a Lie algebra. Let $
G_{+1}$=$\{(^0_0\,^B_0)\},G_{-1}$= $\{(^0_C\,^0_0)\}$. Then $G$ has a $
\Z_2$-consistent \Z\,grading $G$=$G_{-1}\oplus${}$G_0\oplus${}$G_{+1}$,$
G_{\ol0}$=$G_0,$ $G_{\ol1}$=$G_{-1}\oplus${}$G_{+1}$.\vs{-2pt}\par The
Cartan subalgebra $H$ consisting of diagonal $(m$+$n$+2)$\times$($m$+$n
$+2) matrices of zero supertrace has dimension $m$+$n$+1. The weight
space $H$* is the dual of $H$, spanned by the forms $\e_a$ ($a$=1,\dots
,$m$+1), $\d_b$ ($b$=1,\dots,$n$+1), where $\e_a$: $x${}$\rar${}$A_{aa}
$, $\d_b$: $x${}$\rar${}$D_{bb}$ for $x$=($^A_C\,^B_D$), with $\sum_{a=
1}^{m+1}\e_a${}$-${}$\sum_{b=1}^{n+1}\d_b$=0; it has an inner product
derived from the Killing form that $\la\e_a|\e_b\ra$=$\d_{ab}$, $\la\e_
a|\d_b\ra$=0, $\la\d_a|\d_b\ra$=$-${}$\d_{ab}$, where $\d_{ab}$ is the
Kronecker symbol. Let $\D,\D_0,\D_1$ be sets of roots, even, odd roots
respectively, $e(\a)$ the root vector corresponding to $\a$. $G$ has a
root space\vs{-2pt} decomposition $G$=$H${}$\oplus\oplus_{\a\in\D}\C e(
\a)$ with the roots and root vectors given by\vs{-7pt} $$\begin{array}
{rclll}\e_a\!-\!\e_b&\lra&E_{ab}&(1\!\zl\!a,b\!\zl\!m\!+\!1,a\!\zn\!b)&(
\mbox{even}),\\\d_a\!-\!\d_b&\lra&E_{m+a+1,m+b+1}&(1\!\zl\!a,b\!\zl\!n
\!+\!1,a\!\zn\!b)&(\mbox{even}),\\\e_a\!-\!\d_b&\lra&E_{a,m+b+1}&(1\!
\zl\!a\!\zl\!m\!+\!1,1\!\zl\!b\!\zl\!n\!+\!1)&(\mbox{odd}),\\\d_a\!-\!
\e_b&\lra&E_{m+a+1,b}&(1\!\zl\!a\!\zl\!n\!+\!1,1\!\zl\!b\!\zl\!m\!+\!1)&
(\mbox{odd}),\end{array}$$ \vs{-12pt}\nl where\vs{-2pt} $E_{ab}$ is the
matrix with entry 1 at $(a,b)$ and 0 otherwise. We shall find it
convenient to use a notation for roots somewhat different to that in
previous papers.$^{2,15,16}$ Define sets\vs{-6pt} $$\mbox{$I_1$=\{$\ol
m$,${\sc\dots}$,$\ol1$\},\ $I_2$ =\{1,${\sc\dots}$,$n$\},\ $I$=$I_1${}
$\cup$\{0\}$\cup${}$I_2,$ where $\ol i$=$-${}$i,i${}$\zi${}$\Z_+$.}$$\vs
{-12pt}\nl Choose a basis for $H$: $h_i$=$E_{m+i+1,m+i+1}${}$-${}$E_{m+
i+2,m+i+2},i${}$\zi${}$I_1${}$\cup${}$I_2,h_0$=$E_{m+1,m+1}$+$E_{m+2,m+
2}.$ The simple roots in $H$* are: $\a_i$=$\e_{m+i+1}${}$-${}$\e_{m+i+2}
,i${}$\zi${}$I_1,\a_0$=$\e_{m+1}${}$-${}$\d_1,\a_i$=$\d_i${}$-${}$\d_{i+
1},i${}$\zi${}$I_2.$ Thus $\a_0$ is the only odd simple root. The
corresponding Dynkin diagram is\vs{-1pt}\nl\hs{150pt} $\vbox{\hbox{o---
o---$\cdots\cdots$---o---$\otimes$---o---$\cdots$---o---o}\vs{-6pt}
\hbox{$\a_{_{\sc\ol m}}$ $\a_{\ol{\sc m-1}}$\ \ \ \ \ \ $\a_{\ol1}$\ \
$\a_0$\ \ $\a_1$ \ \ \ \ \ $\a_{n-1}$ $\a_n$}}$ \hfill(2.1)\vs{-1pt}\nl
with $I_1,I_2$ corresponding to $sl(m$+1$),sl(n$+1). The symmetric
inner product satisfies\vs{-5pt} $$\matrix{\la\a_i|\a_i\ra\!=\!2,i\!\in
\!I_1,\hfill&\la\a_0|\a_0\ra\!=\!0,\hfill&\la\a_i|\a_i\ra\!=\!-2,i\!\in
\!I_2,\hfill\cr\la\a_{i-1}|\a_i\ra\!=\!-1,i\!\in\!I_1,\hfill&\la\a_0|\a_
{\pm1}\ra\!=\!\pm1,\hfill&\la\a_i|\a_{i+1}\ra\!=\!1,i\!\in\!I_2,\hfill
\cr}\eqno(2.2)$$ \vs{-11pt}\nl and $\la${}$\a_i|\a_j${}$\ra$=$0,j${}$
\ne${}$i,i${}$\pm$1 and $h_i(\a_j)$=$\a_j(h_i)$=$\la${}$\a_i|\a_j${}$
\ra,i${}$\le$0 or $-${}$\la${}$\a_i|\a_j${}$\ra,i${}$>$0. Define\vs
{-10pt} $$\l,\mu\!\in\!H\mbox{*}{\sc\!}:\l\!\ge\!\mu\Lra\l\!-\!\mu\!=\!
\rb{-6pt}{\mbox{$^{\sc\sum}_{i\in I}$}}k_i\a_i\mbox{ with all }k_i\!\ge
\!0,\eqno(2.3)$$\vs{-19pt}\nl a partially order on $H$*. Let $\D^\pm(\D_
0^\pm,\D_1^\pm)$ be sets of positive/negative roots (even, odd roots).
Elements of $\D^+$ are sums of simple roots corresponding to connected
subdiagrams of (2.1). Let $\a_{ij}$=${\sc\sum}^j_{k=i}\a_k,$ then $\D^
\pm_0$=$\{\pm\a_{ij}|i${}$\zl${}$j,i$,$j${}$\zi${}$I_1$ or $i$,$j${}$
\zi${}$I_2$\}, $\D^\pm_1$=$\{\pm\a_{ij}|i${}$\zi${}$I_1${}$\cup$\{0\},$
j${}$\zi$\{0\}$\cup${}$I_2$\}. The root vectors $e_{ij}$=$e(\a_{ij}),f_
{ij}$=$f(\a_{ij})$=$e(-\a_{ij})$ and the elements $h_{ij}$ of $H$ are\vs
{-8pt} $$\mbox{$e_{ij}$=$E_{m+i+1,m+j+2}$, $f_{ij}$=$E_{m+j+2,m+i+1}$,
$h_{ij}$=$E_{m+i+1,m+i+1}${}$-${}$ (-1)^{\si_{ij}}E_{m + j + 2,m + j +
2}$},$$\vs{-18pt}\nl where $\si_{ij}$=0 or $1\Lra\a_{ij}$ is even or
odd. Set $e_i$=$e_{ii},f_i$=$f_{ii}$. The above implies $h_i$=$h_{ii}$
and\vs{-8pt} $$h_{ij}\!=\!\SUM^j_{k=i}h_k,\,i,j\!\zi\!I_1\mbox{ or }i,j
\!\zi\!I_2,\mbox{ \ and \ }h_{ij}\!=\!\SUM^0_{k=i}h_k\!-\!\SUM^j_{k=1}h_
k,\,i\!\le\!0,j\!\ge\!0.$$\vs{-18pt}\nl The set $\{e_{ij},f_{ij},h_i|i,
j\zi I,i\zl j\}$ yields a basis for $G$, with the following nontrivial
relations: $$\matrix{[e_{ij},e_{j+1,\ell}]=e_{i\ell},\ \ \ [f_{ij},f_{j+
1,\ell}]=-f_{i\ell},\ \ \ [e_{ij},f_{ij}]=h_{ij},\hfill\vs{4pt}\cr[e_
{ij},f_{ik}]=\left\{^{\dis-(-1)^{\si_{ij}\si_{ik}}f_{j+1,k}\mbox{ if }j<
k,}_{\dis-(-1)^{\si_{ij}\si_{ik}}e_{k+1,j}\mbox{ if }j>k,}\right.\ \ \ [
e_{ik},f_{jk}]=\left\{^{^{\dis e_{i,j-1}\mbox{ if }i<j,}}_{_{\dis f_{j,
i-1}\mbox{ if }i>j,}}\right.\hfill\vs{4pt}\cr[h_{ij},e_{k\ell}]\!=\!\mu
e_{k\ell},[h_{ij},f_{k\ell}]\!=\!-\mu f_{k\ell},\mu\!=\!\d_{i,k}\!-\!\d_
{i,\ell+1}\!-\!(-1)^{\si_{ij}}\d_{j,k-1}\!+\!(-1)^{\si_{ij}}\d_{j,\ell}
\,.\hfill\cr}\eqno(2.4)$$ Set $G^\pm_0=\mbox{span}\{e(\a)|\a\zi\D^\pm_0
\},G\pm_1=\mbox{span}\{e(\b)|\b\zi\D^\pm_1\},G^\pm=G^\pm_0\oplus G^\pm_
1.$ Note that $G^\pm_1=G_{\pm1},G_{\ol0}=G^-_0\oplus H\oplus G^+_0,G=G^-
\oplus H\oplus G^+.$ Let ${\bf U}(G)$ be the universal enveloping
algebra of $G$, ${\bf U}(G')$ that of its subalgebras $G'$ which is $H
$-diagonalizable. Denote by ${\bf U}(G')_\eta$ the subspace of weight $
\eta$. The PBW theorem can be extended to Lie superalgebras:$^{4,7}$
\par\ni{\bf Theorem 2.1.} Let $y_1,{\sc\cdots},y_{\ssc M}$ be a basis
of $G_{\ol0}$ and $z_1,{\sc\cdots},z_{\ssc N}$ be that of $G_{\ol1}$.
The elements of the form $(y_1)^{k_1}{\sc\cdots}(y_{\ssc M})^{k_M}z_{i_
1}{\sc\cdots}z_{i_s}$, where $k_i\zg 0$ and $1\zl i_1<{\sc\cdots}<i_s
\zl N$, form a basis of ${\bf U}(G)$. \qed\par\ni For $\l\zi H$*,
define its {\sl Dynkin labels} to be $a_i=\l(h_i),i\zi I$. These
uniquely determine $\l$, which can then be represented as $\l=[a_{\ol m}
,{\sc\cdots},a_{\ol1};a_0;a_1,{\sc\cdots},a_n]$. $\l$ is called {\sl
dominant} if $a_i\zg 0$ for all $i\zn 0$, {\sl integral} if $a_i\zi\Z$
for all $i\zn 0$. The following convention will be useful later. \par\ni
{\sl Convention 2.2.} If $\G$ denotes any quantity relating\vs{-2pt} to
$G=sl(m+1/n+1)$, then $\G^{(m'/n')}$ denotes the same quantity relating
to $sl(m'+1/n'+1)$. Thus $\G^{(m/n)}=\G$.\qed\es {\bf III. THE
KAC-MODULES}\vs{-1pt}\par\ni \vs{-1pt}Let $V^0(\L)$ be the simple $G_{
\ol0}$-module with integral dominant highest weight $\L$ and vector $v_
\L$. Extend\vs{-1pt} $V^0(\L)$ to be a $G_{\ol0}\oplus G_{+1}$ module
by setting $G_{+1}V^0(\L)=0$. The {\sl Kac-module}$^{\,6}$ is\vs{-7pt}
$$\VBL=\mbox{Ind}^G_{G_0\oplus G_{+1}}V^0(\L)={\bf U}(G)\otimes_{G_0
\oplus G_{+1}}V^0(\L).$$ \vs{-15pt}\nl By Theorem 2.1,\vs{-1pt} ${\bf U}
(G)={\bf U}(G_{-1})\otimes{\bf U}(G_0)\otimes{\bf U}(G_{+1})$. It
implies $\VBL\cong{\bf U}(G_{-1})\otimes V^0(\L).$ We summarize some
well known properties of $\VBL$; more details can be found in {\sl
Refs.} 2, 6. By definition, it is a $2^{(m+1)(n+1)}\mbox{dim}V^0(\L)$
dimensional highest weight module generated by the highest weight
vector $v_\L$, indecomposable and $H$-diagonalizable and it contains a
maximal submodule $M=\{v\zi\VBL|v_\L\not\in{\bf U}(G)v\}$, such that $V(
\L)=\VBL/M$ is a finite-dimensional simple module with highest weight $
\L$. Define $\rho=\rho_0-\rho_1,\rho_0={1\over2}{\sc\sum}_{\a\in\D^+_0}
\a,\rho_1={1\over2}{\sc\sum}_{\b\in\D^+_1}\b.$\par\ni {\sl Definition
3.2.} $\L,\VBL,V(\L)$ are called {\sl typical} if $\la\L+\rho|\b\ra\ne
0$ for all $\b\zi\D^+_1$. If $\b\zi\D_1^+$ such that $\la\L+\rho|\b\ra=
0$,\,then $\L,\VBL,V(\L)$ are called {\sl atypical} and $\b$ is an {\sl
atypical root} for $\L$. If there exist precisely $r$ distinct atypical
roots for $\L$, we call $\L,\VBL,V(\L)$ {\sl $r$-fold atypical}.\qed
\par\ni{\bf Theorem 3.3 (Kac}$^{\,6}${\bf).} Every finite-dimensional
simple $G$-module is isomorphic to a $V(\L)$, characterized by its
integral dominant highest weight $\L$. $\VBL$ is simple $\Lra$ $\L$ is
typical. \qed\par\ni A {\sl composition series} of $\VBL$ is a sequence
$\VBL=V_0\supset V_1\supset\cdots$ with each $V_i/V_{i+1}$ isomorphic
to some simple module $V(\Si)$, called a {\sl composition factor} of $
\VBL$. A conjecture was made in {\sl Ref.} 2, giving all the
composition factors of $\VBL$. We aim to prove the existence of some of
these composition factors; for this, important concepts are those
defined as follows. \par\ni{\sl Definition 3.4.} A vector $v$ in a $G
$-module $V$ is called {\sl weakly $G$-primitive} if there exists a $G
$-submodule $U$ of $V$ such that $v\not\in U$ and $G^+v\subset U$. If $
U=0$, $v$ is called {\sl $G$-primitive}.\qed\par\ni We are only
concerned with finite-dimensional modules. Thus, weakly $G_{\ol0}
$-primitive vectors is in fact $G_{\ol0}$-primitive and integral
dominant. A {\sl cyclic module} is an indecomposable module generated
by a weakly primitive vector. A weakly primitive vector $v$ will
determine a cyclic submodule\,${\bf U}(G)v$\,and\,a\,composition
factor. An\,important\,construct\,in\,classifying\,composition factors
is the atypicality matrix $A(\L)$.$^{2,15,16}$ First, introduce the
shorthand notation:\vs{-7pt} $$\b_{bc}=\e_b-\d_c=\a_{_{\sc\ol{m-b+1},c-
1}},\ 1\zl b\zl m+1,1\zl c\zl n+1.$$ \vs{-13pt}\nl {\sl Definition
3.5.} The\,atypicality\,matrix$\,A(\L)\,$is\,the\,($m$+1)$\times$($n
$+1)$\,$matrix with$\,$($b$,$c$)-entry$\,A(\L)_{bc}$ =$\la\L$+$\rho|\b_
{bc}\ra$=$\SUM^0_{_{\sc k=\ol{m-b+1}}}a_k${}$-${}$\SUM^{c-1}_{k=1}a_k$+$
m${}$-${}$b${}$-${}$c$+2. For example, for $G$=$sl$(4/5$),\L
$=[100;0;1000],\vs{-4pt} \nl\hs{130pt}$A(\L)={\sc\pmatrix{\sc 4\ 2\ 1\
0\ \ol1\vs{-5pt}\cr\sc2\ 0\ \ol1\ \ol2\ \ol3\vs{-5pt}\cr\sc1\ \ol1\
\ol2\ \ol3\ \ol4\vs{-5pt}\cr\sc0\ \ol2\ \ol3\ \ol4\ \ol5\cr}}.$\qed\vs
{2pt}\nl\ni Inspection\,of\,$A({\sc\!}\L{\sc\!})
$\,tells\,immediately\,whether\,or\,not\, ${\!}\L$
is\,atypical\,and\,which\,are\,the\,atypical\,roots
since\,they\,correspond\,to\,zero\,entries\,of $A({\sc\!}\L{\sc\!})
$.\,In\, ${\sc\!}$above,\,$\L{\sc\!}$
is\,3-fold\,atypical\,with\,atypical\,roots\,$\b_{41}$, $\b_{22}$,$\b_
{14}$.\,The\,properties\,of\,$A({\sc\!}\L{\sc\!})
$\,have\,been\,studied\,in\,detail\,in\,{\sl Ref.} 2.
We\,summarize\,some\,here. \vs{-2pt}\par\ni{\sl Lemma 3.6.} (i) Let $\L=
[a_{\ol m},\cdots,a_{\ol1};a_0;a_1,\cdots,a_n]$; then\vs{-7pt} $$\begin
{array}{rcl}A(\L)_{bc}-A(\L)_{b+1,c}&=&a_{_{\sc\ol{m-b+1}}}+1,\ 1\zl b
\zl m,1\zl c\zl n+1,\\ A(\L)_{m+1,1}&=&a_0,\\ A(\L)_{bc}-A(\L)_{b,c+1}&=
&a_c+1,\ 1\zl b\zl m+1,1\zl c\zl n.\end{array}\eqno(3.1{\rm a})$$ \vs
{-13pt}\nl\hs{3ex}(ii) An atypicality matrix $A(\L)$ satisfies $A(\L)_
{bc}+A(\L)_{de}=A(\L)_{be}+A(\L)_{dc}.$ {\it Vice versa}, any $(m+1)
\times(n+1)$ matrix satisfying this condition for all pairs $(b,c),\,(d,
e)$ with $1\zl b,d\zl m+1$ and $1\zl c,e\zl n+1$ is the atypicality
matrix of a unique element $\L\zi H$*.\par\ni \hs{3ex}(iii) $\L$ is
dominant $\Lra$\vs{-5pt} $$\matrix{A(\L)_{bc}-A(\L)_{b+1,c}-1\zg0,&1\zl
b\zl m,1\zl c\zl n+1,\cr A(\L)_{bc}-A(\L)_{b,c+1}-1\zg0,&1\zl b\zl m+1,
1\zl c\zl n.\cr}\eqno(3.1{\rm b})$$ Moreover, $\L$ is integral if the
expressions on the {\sl l.h.s.} of (3.1b) are all integers.\qed\par\ni
For atypical modules, the highest weight $\L$ must be integral dominant
and $a_0$ is an integer since at least one of the entries of $A(\L)$ is
zero. Lemma 3.6 implies that the zeros of $A(\L)$ lie in distinct rows
and columns, and that one zero lies to the right of another $\Lra$ it
lies above it. Thus the atypical roots are commensurate with respect to
ordering (2.3). If $\L$ is $r$-fold atypical, we label the atypical
roots $\g_1<{\sc\cdots}<\g_r$. It follows that if $1\zl s,t\zl r$ and $
x_{st}$ is the entry in $A(\L)$ at the intersection of the column
containing the $\g_s$ zero with the row containing the $\g_t$ zero,
then $x_{st}\zi\Z_+\bs\{0\}$ for $s<t$ and $x_{ts}=-x_{st}$. Therefore $
A(\L)$ has the form:\nl\hs{120pt} $A(\L)=\pmatrix{_{\cdots}{}^{\ \ssc
\vdots}_{x_{1r}}{}_{\cdots}{}^{\ \ssc\vdots}_{x_{2r}}{}_{\cdots}{}^{\
\ssc\vdots}_{x_{3r}}{}_{\cdots}{}^{\ \ssc\vdots}_{\ 0\ }{}_{\cdots}\vs{-
5pt}\cr_{\cdots}{}^{\ \ssc\vdots}_{x_{13}}{}_{\cdots}{}^{\ \ssc\vdots}_
{x_{23}}{}_{\cdots}{}^{\ \ssc\vdots}_{\ 0\ }{}_{\cdots}{}^{\ \ssc\vdots}
_{\ol{x_{3r}}}{}_{\cdots}\vs{-5pt}\cr_{\cdots}{}^{\ \ssc\vdots}_{x_{12}}
{}_{\cdots}{}^{\ \ssc\vdots}_{\ 0\ }{}_{\cdots}{}^{\ \ssc\vdots}_{\ol
{x_{23}}}{}_{\cdots}{}^{\ \ssc\vdots}_{\ol{x_{2r}}}{}_{\cdots}\vs{-5pt}
\cr_{\cdots}{}^{\ \ssc\vdots}_{\ 0\ }{}_{\cdots}{}^{\ \ssc\vdots}_{\ol
{x_{12}}}{}_{\cdots}{}^{\ \ssc\vdots}_{\ol{x_{13}}}{}_{\cdots}{}^{\
\ssc\vdots}_{\ol{x_{1r}}}{}_{\cdots}\cr}.$\hfill(3.2)\par\ni Denote $h_
{st}$ the hook length between the zeros corresponding to $\g_s,\g_t$,
i.e., the number of steps to go from the $\g_s$ zero via $x_{st}$ to
the $\g_t$ zero with the zeros themselves included in the count. An
important concept in the classification of composition factors is the
following.$^2$ \par\ni{\sl Definition 3.7.} Let $\L$ be $r$-fold
atypical with atypical roots $\{\g_1,{\sc\cdots},\g_r\}$. For $1\zl s<t
\zl r$: \nl\hs{3ex}\rb{13pt}{\mbox{$\matrix{ \mbox{(i) $\g_s,\g_t$ are
{\sl normally related} ($n$)}\hfill& \Lra&x_{st}>h_{st}-1;\hfill\cr
\mbox{(ii) $\g_s,\g_t$ are {\sl quasi-critically related} ($q$)}\hfill&
\Lra&x_{st}=h_{st}-1;\hfill\cr \mbox{(iii) $\g_s,\g_t$ are {\sl
critically related} ($c$)}\hfill& \Lra&x_{st}<h_{st}-1.\hfill\cr}$}}\qed
\par\ni It is straightforward to show that the $q$-relation is
transitive, i.e., if $\g_s,\g_t$ are $q$-related and $\g_t,\g_u$ are $q
$-related, then $\g_s,\g_u$ are $q$-related. \par\ni{\sl Definition
3.8.} The {\sl $nqc$-type} ({\sl atypicality type}) of an $r$-fold
atypical $\L$ is a triangular array\nl\hs{100pt} $nqc(\L)=\ ^{^{\dis s_
{\sc1r}\,{\sc\cdots}\,{\dis s}_{\sc sr}\,{\sc\cdots}\,{\dis s}_{\sc tr}
\,{\sc\cdots}\,\,0}_{{_{{\dis s}_{\sc1t}}^{\,\vdots}}\,_{\cdots}\,{_{{
\dis s}_{\sc st}}^{\,\vdots}}\,_{\cdots}\,\,{_{\sc0}^{\vdots\ \ ^{\sc
\cdot^{\sc\cdot^{\sc\cdot}}}}}}}_{^{{_{{\dis s}_{\sc1s}}^{\,\vdots}}\,_
{\cdots}\,\,{_{\dis0}^{\vdots\ \ ^{\sc\cdot^{\sc\cdot^{\sc\cdot}}}}}}_{
\,{_{\dis0}^{\vdots\ \ ^{\sc\cdot^{\sc\cdot^{\sc\cdot}}}}}}}\ \ ,\ \ \
s_{st}\in\{n,q,c\},$\nl where the zeros correspond to $\{\g_1,\cdots,\g_
r\}$ and $s_{st}=n,q,c\Lra\g_s,\g_t$ are $n$-, $q$-, $c$-related.\qed
\par\ni It was conjectured$^{\,2}$ that the number and nature of
composition factors of $\VBL$ depends only on $nqc$-type of $\L$; if a
weight $\L$ of $sl(m+1/n+1)$ and a weight $\L'$ of $sl(m'+1/n'+1)$ have
the same $nqc$-type, then there is a 1 - 1 correspondence between the
composition factors of $\VBL$ and $\ol{V}(\L')$. More precisely it was
conjectured that the composition factors of $\VBL$ are in 1 - 1
correspondence with certain codes $\Si^c$ which are determined from $
nqc(\L)$, and which in turn determine the highest weights $\Si$ of the
corresponding composition factors $V(\Si)$. \par\ni{\sl Definition
3.9.} A {\sl code} $\Si^c$ for $\L$ is an array of length $r$, each
element of the array consisting of a non-empty column of increasing
labels taken from $\{0,{\sc\cdots},r\}$. The 1st element of a column is
called the {\sl top label}. $\Si^c$ must satisfy the rules: \nl\hs
{3ex}(i) The top label of column\,$s$\,can be\,$0,s$ or $a$\,with $s${}$
<${}$a$; the 1st case can occur only if column $s$ is zero, while the
last case can occur only if $nqc(\L)_{st}$=$q$ with $a$ the top label
of column $t$. \nl\hs{3ex}(ii) Let $s<t,nqc(\L)_{st}={\sc\cdots}=nqc(\L)
_{t-1,t}=c$. If the top label of column $t$ is $a:\,t\zl a$, then $a$
must appear somewhere below the top entry of column $s$. \nl\hs
{3ex}(iii) If $s$ appears in any column then the only labels which can
appear below $s$ in the same column are those $t:\,s<t$, for which $t$
is the top label of column $t$ and $nqc(\L)_{st}=c$. \nl\hs{3ex}(iv) If
the label $s$ appears in more than one column and $t$ appears
immediately below $s$ in one such column, then it must do so in all
columns containing $s$. \nl\hs{3ex}(v) Let $s<t<u,\,nqc(\L)_{st}=q,nqc(
\L)_{tu}=q$ (so, $nqc(\L)_{su}=q$). If the top label of column $s$ is
the same as that of column $u$ and it is non-zero then the top label of
column $t$ is not 0. \nl\hs{3ex}(vi) Let $s<t<u<v$ with top labels $a,b,
a,b$ respectively, $a\zn0\ne b$. If $a<b$ then columns $s$ and $u$ must
contain $b$; if $a>b$ then columns $t$ and $v$ must contain $a$.\qed\par
\ni As an example, consider $\L=[00020;0;0210]$ for $sl(6/5)$. A
straightforward computation gives\nl\hs{120pt} $A(\L)={\ssc\pmatrix{\sc
7\ 6\ 3\ 1\ 0\vs{-5pt}\cr\sc 6\ 5\ 2\ 0\ \ol1\vs{-5pt}\cr\sc 5\ 4\ 1\
\ol1\ \ol2\vs{-5pt}\cr\sc 4\ 3\ 0\ \ol2\ \ol3\vs{-5pt}\cr\sc 1\ 0\ \ol3
\ \ol5\ \ol6\vs{-5pt}\cr\sc 0\ \ol1\ \ol4\ \ol6\ \ol7\vs{-5pt}\cr}}\ \ ,
\ \ \ \ \ nqc(\L)=\ ^{^{\sc c\ c\ c\ c\ 0}_{\sc c\ q\ c\ 0}}_{^{^{\sc q
\ n\ 0}_{\sc c\ 0}}_{\sc 0}}\ \ \ .$\hfill(3.3)\vs{2pt}\nl $\L$ is
5-fold atypical with atypical roots $\g_1=\b_{61},\g_2=\b_{52},\g_3=\b_
{43},\g_4=\b_{24},\g_5=\b_{15}$. Using the rules in Definition 3.9, we
find the following 15 codes $\Si^c$: $$\matrix{^{0\ 0\ 0\ 0\ 0}&^{1\ 0\
0\ 0\ 0}&^{\sc1\ 2\ 0\ 0\ 0}_{\sc2}&^{0\ 0\ 3\ 0\ 0}&^{1\ 0\ 3\ 0\ 0}&^
{1\ 2\ 3\ 0\ 0}_2&^{3\ 0\ 3\ 0\ 0}\hfill&\!\!\!\!\!^{0\ 0\ }{}^3_4{}^{\
4\ 0}\hfill\vs{5pt}\cr^{1\ 0\ }{}^3_4{}^{\ 4\ 0}&^{1\ 2\ }_2{}^{3\ 4\ 0}
_4&^{1\ 4\ }_4{}^{3\ 4\ 0}_4&^{3\ 4\ }_4{}^{3\ 4\ 0}_4&^{1\ 2\ }_{^{
\sc2\ 5}_{\sc5}}{}^{3\ 4\ 5}_{^{\sc4\ 5}_{\sc5}}&^{1\ 4\ }_{^{\sc4\ 5}_
{\sc5}}{}^{3\ 4\ 5}_{^{\sc4\ 5}_{\sc5}}&^{3\ 4\ }_{^{\sc4\ 5}_{\sc5}}{}^
{3\ 4\ 5}_{^{\sc4\ 5}_{\sc5}}\ \ .\cr}\eqno(3.4)$$ For $1\zl s\zl r$,
the $s$-th column of a code corresponds to the $s$-th atypical root $\g_
s$. From definition, we see that if $\g_s,\g_t$ are $q$-related and the
top entry $a$ of column $s$ is non-zero, then $a$ may also be the top
entry of column $t$. In such a case, we say that $\g_s,\g_t$ are {\sl
linked}. In the example, $\g_1,\g_3$ are $q$-related, and they are
linked in code $(3 0 3 0 0)$, whereas in code $(1 0 3 0 0)$, they are
not. Thus, where $\g_s,\g_t$ are $q$-related, there will be codes in
which they are linked, and codes in which they are not. This leads the
following definition.\par\ni {\sl Definition 3.10.} A code $\Si^c$ is a
{\sl linked code} if there exist $\g_s,\g_t$ which are linked, i.e.,
columns $s$ and $t$ have the same non-zero top entry. Otherwise, it is
called an {\sl unlinked code}.\qed\par\ni It follows from rule (i) that
if $nqc(\L)$ contains no $q$, then all codes are unlinked. Next, we see
from rule (ii) that if, for $s<t,$ $\g_s,{\sc\cdots},\g_{t-1}$ are $c
$-related to $\g_t$ and if the top label $a$ of column $t$ of a code is
non-zero, then $a$ must appear somewhere below the top entry of the $s
$-th column, and so also of the $(s$+1)-th,${\sc\cdots}$,$(t$-1)-th
column, of that code; we say\,$\g_t$\,{\sl wraps} $\g_s$.\,Unlike
links, wraps must be made. In the example, $\g_1,\g_2$ are $c$-related,
and each code in which the 2nd column is non-zero, e.g., $^{1\ 2\ 0\ 0\
0}_{2},\ ^{1\ 4\ }_{4}{}^{3\ 4\ 0}_{4}$, the top label of the 2nd
column occurs below the top entry in the 1st column; i.e., $\g_2$ wraps
$\g_1$. Similarly, $\g_1,{\sc\cdots},\g_4$ are $c$-related to $\g_5$,
and in each code with non-zero top label in the 5th column, that label
occurs \vs{-3pt}below the top entry in each of the first 4 columns, as
in \ \rb{2pt}{$^{3\ 4\ }_{^{\sc4\ 5}_{\sc5}}{}^{3\ 4\ 5}_{^{\sc4\ 5}_{
\sc5}}$}. \vs{-4pt}It may happen, on the other hand, that for $s<t,\g_s,
\g_t$ are $c$-related but $\g_u,\g_t$ are not $c$-related for some $u,
\,s<u<t$, as in (3.3) $\g_1,\g_4$ are $c$-related but $\g_2,\g_4$ are $
q$-related. Correspondingly, in code $^{1\ 2\ }_2{}^{3\ 4\ 0}_4$, $\g_4
$ does not wrap $\g_1$. However, in code $^{1\ 4\ }_4{}^{3\ 4\ 0}_4$, $
\g_4$ does appear to wrap $\g_1$. This is because in this code, as
opposed to the previous one, $\g_4$ is linked to $\g_2$, so the top
label of the 2nd column is the same (i.e., 4) as that of the 4th
column, and therefore, since $\g_2$ must wrap $\g_1$ in all codes in
which the 2nd column is non-zero, this entry must appear below the top
entry in the 1st column. Thus $\g_4$ wraps $\g_1$ only because of the
presence of an intermediate link; we shall use the term {\sl link wrap}
rather than wrap to describe this. From the discussion we see that, in
general, the presence of a $q$ rather than an $n$ in $nqc(\L)$
increases, whereas the presence of a $c$ rather than an $n$ decreases,
the number of codes for $\L$. Thus, for $r=2$, $\L$ has $3,4,5$ codes $
\Lra nqc(\L)=\,^{c\ 0}_0,\ ^{n\ 0}_0,\ ^{q\ 0}_0$. \par\ni {\sl
Definition 3.11.} In a code $\Si^c$ for $\L$, we say that $\g_s$ is {\sl
connected} to $\g_t$ and write \lk{\g_s}{\g_t} if the $s$-th and $t$-th
columns of $\Si^c$ contain a common non-zero entry.\qed\par\ni Thus for
the code $\Si^c=\,^{1\ 2\ 3}_{3\ 3}$, \lk{\lk{\g_1}{\g_2}}{\g_3},
whereas for the code $\Si^c=1\,2\,0$, $\g_1\mbox{\rb{5pt}{$_{^{\dis\not
\sim}_{\ssc\Si}}$}}\g_2$. Clearly, \lk{}{} is an equivalent relation on
$\{\g_s\,|$ the $s$-th column of $\Si$ is non-zero$\}$.\par\ni {\sl
Definition 3.12.} A code $\Si^c$ is called {\sl indecomposable} if $\{
\g_s\,|$ the $s$-th column of $\Si^c$ is non-zero$\}$ is an equivalent
class for the relation \lk{}{}; otherwise $\Si^c$ is called {\sl
decomposable}.\qed\par\ni For example, $^{1\ 2\ 3\ 0}_{3\ 3}$ is
indecomposable, whereas $^{1\ 2\ }_2{}^{3\ 4}_4$ is decomposable. If $
\Si^c$ is decomposable, then we can write $\Si^c=\Si_1^c\Si_2^c\cdots
\Si_s^c\cdots$, where each $0\cdots0\Si_s^c0\cdots0$ (with 0's in the
appropriate positions) is indecomposable. For instance, $\Si^c=\,^{1\ 2
\ }_2{}^{3\ 4}_4$ can be written as $\Si^c=\Si_1^c\Si_2^c$ where $\Si_1^
c=\,^{1\ 2}_2,\,\Si_2^c=\,^{3\ 4}_4$, with $^{1\ 2\ 0\ 0}_2$ and $^{0\
0\ }{}^{3\ 4}_4$ being themselves indecomposable unlinked codes for $\L
$. Without confusion we will simply denote $0\cdots0\Si_s^c0\cdots0$ by
$\Si_s^c$.\es {\bf IV. SOUTH WEST CHAINS OF $A(\L)$}\par\ni To obtain
the highest weights of those composition factors of $\VBL$
corresponding to unlinked codes, we need to construct south west
chains.$^2$ Denote by $$D=\{(b,c)|1\zl b\zl m+1,1\zl c\zl n+1\},\ \ \ \
\G_\L=\{(b,c)|A(\L)_{bc}=0\}.\eqno(4.1)$$ the set of, respectively,
positions, positions of zeros, of $A(\L)$ and define $\wh K=\{\b_{bc}\,|
\,(b,c)\in K\}$ for any subset $K$ of $D$. In particular, $\wh D=\D_1^+
$ and $\wh\G_\L=\{\g_1,\ldots,\g_r\}$.\par\ni {\sl Definition 4.1.} (i)
For 1$\zl${}$s${}$\zl${}$r$, let ($b_s$,$c_s$)$\in${}$\G_\L$ be the
position corresponding to $\g_s$. The {\sl extended west chain} $W^e_\L(
s)$ emanating from ($b_s$,$c_s)$ is a sequence of positions in $D$
starting at ($b_s$,$c_s)$ and extending in a westerly or south-westerly
direction until it reaches the 1st column or it cannot extend further
without leaving $A(\L)$ by passing below its bottom row. For all $t$
with 1$\zl${}$t${}$\zl${}$c_s,$ $W^e_\L(s)$ has exactly one element in
the $t$-th column provided that the row of this element lies within $A(
\L)$. For 1$\zl${}$t${}$\zl${}$c_s${}$-$1, the row of the position in
the $t$-th column is $a_t$ rows below the row of the position in the $(
t+1)$-th column, where $a_t$ is a Dynkin label of $\L$; if this is not
possible, i.e., if this row would be the $M$-th row where $M${}$>${}$m
$+1, then $W^e_\L(s)$ ends in the $t$-th column, i.e., has no position
to the left of the $t$-th column. Thus $W^e_\L(s)$ is the set\vs{-9pt}
$$W^e_\L(s)=D\cap\{(b,c)|1\le c\le c_s,b=b_s+\rb{4pt}{\mbox{$^{c_{_{\sc
s}}-1}_{^{\ \sum}_{\sc\ t=c}}$}}a_t\}.\eqno(4.2{\rm a})$$ \vs{-17pt}\nl
\hs{3ex}(ii) Similarly, the {\sl extended south chain} $S^e_\L(s)$
emanating from $(b_s,c_s)$ is the set\vs{-9pt} $$S^e_\L(s)=D\cap\{(b,c)|
b_s\le b\le m+1,c=c_s-\rb{4pt}{\mbox{$^{\,b-1}_{^{\,\sum}_{\sc t=b_{\sc
s}}}$}}a_{_{\sc\ol{m-t+1}}}\}.\eqno(4.2{\rm b})$$ \vs{-17pt}\nl \hs
{3ex}(iii) The {\sl extended south west chain} emanating from ($b_s$,$c_
s$) is $SW^e_\L$($s$)=$W^e_\L$($s$)$\cup${}$S^e_\L$($s$).\qed\par \ni
{\sl Definition 4.2.} (i) For $1\le s\le r$, the {\sl south west chain}
$SW_\L(s)$ emanating from $(b_s,c_s)$ consists of all positions of $SW^
e_\L(s)$ which are above and to the right of {\bf all} points of
intersection of the chains $W^e_\L(s)$ and $S^e_\L(s)$ of $A(\L)$, with
the additional requirement that if $S^e_\L(s)$ starts off at $(b_s,c_s)
$ by immediately going above $W_\L^e(s)$, then $SW_\L(s)$ consists
solely of $(b_s,c_s)$. The {\sl west}, {\sl south subchain} of $SW_\L$($
s$) are, respectively, $W_\L$($s$)=$W_\L^e$($s$)$\cap${}$SW_\L$($s$), $
S_\L$($s$)=$S_\L^e$($s$)$\cap${}$SW_\L$($s$). \nl\hs{3ex}(ii) $SW_\L=
\cup^r_{s=1}SW_\L(s)$ is called the set of all south west chains.\qed
\par\ni The construction of chains is facilitated by placing the Dynkin
labels $a_{\ol m},{\sc\ldots},a_{\ol1}$ to the left of the 1st column,
and in between the rows, of $A(\L)$, likewise, $a_1,{\sc\ldots},a_n$
are placed below the bottom row, and in between the columns, of $A(\L)
$. We illustrate chains with 3 examples of doubly atypical $\L$ for $sl(
5/6)$.$ $ In each case we first give $SW_\L^e(s)$, denoting $W_\L^e(s)$
by broken lines, $S^e_\L(s)$ by unbroken lines, and then give $SW_\L(s)
$, denoting the chains by arrows.\par\ni {\sl Example 4.3.} $\L=[1\,1
\,1\,1;\ol 1;1\,0\,0\,1\,0]$. \par\ \ \ \ \ \ \ \ \ \ $A(\L)$ =$\ \ $
\put(0,30){$\sc1$}\put(0,10){$\sc1$}\put(0,-10){$\sc1$}\put(0,-30){$
\sc1$} \put(15,0){\oval(10,76)[l]} \put(20,40){$\sc7$}\put(20,20){$\sc5
$}\put(20,0){$\sc3$}\put(20,-20){$\sc1$} \put(20,-40){$\sc\ol1$}\put
(35,40){$\sc5$} \put(35,20) {$\sc 3$}\put(35,0) {$\sc1$}\put(35,-20){$
\sc\ol1$}\put(35,-40){$\sc\ol3$} \put(25,7){$\sc.$}\put(29,12){$\sc.
$}\put(33,17){$\sc.$} \put(25,-13){$\sc.$}\put(29,-8){$\sc.$}\put
(33,-3){$\sc.$} \put(25,-33){\line(3,4){8}}\put(50,40){$\sc4$} \put
(50,20){$\sc2$}\put(50,0){$\sc0$}\put(50,-20){$\sc\ol2$}\put(50,-40) {$
\sc\ol4$}\put(40,22){$\sc.$}\put(43,22){$\sc.$}\put(46,22){$\sc.$} \put
(40,2){$\sc.$}\put(43,2){$\sc.$}\put(46,2){$\sc.$} \put(40,-13){\line
(3,4){8}}\put(40,-33){\line(3,4){8}} \put(65,40) {$\sc 3$} \put
(65,20){$\sc1$}\put(65,0){$\sc\ol1$} \put(65,-20) {$\sc\ol 3$} \put
(65,-40){$\sc\ol5$}\put(55,22){$\sc.$}\put(58,22){$\sc.$}\put(61,22) {$
\sc.$}\put(55,-13){\line(3,4){8}}\put(80,40){$\sc1$} \put(80,20){$\sc
\ol1$}\put(80,0){$\sc\ol3$}\put(80,-20){$\sc\ol5$} \put(80,-40){$\sc
\ol7$}\put(70,7){\line(3,4){8}} \put(70,27){$\sc.$}\put(74,32){$\sc.
$}\put(79,37){$\sc.$}\put(95,40){$\sc0$} \put(95,20){$\sc\ol2$}\put
(95,0){$\sc\ol4$}\put(95,-20){$\sc\ol6$} \put(95,-40){$\sc\ol8$}\put
(85,42){$\sc.$}\put(88,42){$\sc.$}\put(91,42) {$\sc.$}\put(85,27){\line
(3,4){8}}\put(102,0){\oval(10,76)[r]} \put(28,-55){$\sc1$}\put(43,-55){$
\sc0$}\put(58,-55){$\sc0$}\put(73,-55) {$\sc1$}\put(88,-55){$\sc0$}\put
(115,0){=$\ \ $ \put(0,30){$\sc1$}\put(0,10){$\sc1$}\put(0,-10){$\sc1
$}\put(0,-30){$\sc1$} \put(15,0){\oval(10,76)[l]}\put(20,40){$\sc7$}
\put(20,20){$\sc5$} \put(20,0){$\sc 3$}\put(20,-20){$\sc1$}\put
(20,-40) {$\sc\ol1$}\put(35,40){$\sc5$} \put(35,20){$\sc 3$}\put
(35,0){$\sc1$} \put(35,-20){$\sc\ol1$}\put(35,-40){$\sc\ol3$}\put
(33,17){\vector(-3,-4) {8}}\put(33,-3){\vector(-3,-4){8}}\put
(33,-23){\vector(-3,-4){8}} \put(50,40){$\sc4$}\put(50,20){$\sc2$}\put
(50,0){$\sc0$}\put(50,-20) {$\sc\ol2$}\put(50,-40){$\sc\ol4$}\put
(46,22){\vector(-1,0){6}} \put(46,2){\vector(-1,0){6}}\put
(48,-3){\vector(-3,-4){8}} \put(48,-23){\vector(-3,-4){8}} \put(65,40){$
\sc 3$} \put(65,20){$\sc1$}\put(65,0){$\sc\ol1$}\put(65,-20){$\sc\ol3
$} \put(65,-40){$\sc\ol5$}\put(61,22){\vector(-1,0){6}} \put
(63,-3){\vector(-3,-4){8}}\put(80,40){$\sc1$} \put(80,20){$\sc\ol1$}\put
(80,0){$\sc\ol3$}\put(80,-20){$\sc\ol5$} \put(80,-40){$\sc\ol7$}\put
(78,17){\vector(-3,-4){8}} \put(79,37){\vector(-3,-4){8}}\put(95,40){$
\sc0$} \put(95,20){$\sc\ol2$}\put(95,0){$\sc\ol4$}\put(95,-20){$\sc\ol6
$} \put(95,-40){$\sc\ol8$}\put(91,42){\vector(-1,0){6}} \put
(93,37){\vector(-3,-4){8}}\put(102,0){\oval(10,76)[r]}\put(112,0){.}
\put(28,-55){$\sc1$}\put(43,-55){$\sc0$}\put(58,-55){$\sc0$}\put
(73,-55) {$\sc1$}\put(88,-55){$\sc0$}}\par\ni $SW^e_\L(1)=SW_\L(1)=\{(3,
2),(3,3),(4,1),(4,2),(5,1)\}$; $SW^e_\L(2)=SW_\L(2)=\{(1,5),$ $(1,6),$ $
(2,2),$ $(2,3),$ $(2,4),$ $(2,5),$ $(3,1),$ $(3,4),$ $(4,3),$ $(5,2)\}
$. Here $\g_1,\g_2$ are $c$-related. \qed\vs{-2pt}\par\ni {\sl Example
4.4.} $\L=[2\ 3\ 0\ 2;\ol1;0\ 1\ 1\ 2\ 1]$.\nl\hs{20pt} \ \ \ \ \ \ \ \
\ \ $A(\L)$ =$\ \ $\put(0,30){$\sc2$} \put(0,10){$\sc3$}\put(0,-10){$
\sc0$}\put(0,-30){$\sc2$} \put(15,0){\oval(10,76)[l]}\put(20,40){$\sc10
$} \put(20,20){$\sc7$}\put(20,0){$\sc3$}\put(20,-20){$\sc2$} \put
(20,-40){$\sc\ol1$}\put(35,40){$\sc9$}\put(35,20){$\sc6$}\put(35,0) {$
\sc2$}\put(35,-20){$\sc1$}\put(35,-40){$\sc\ol2$}\put(50,40){$\sc7$}
\put(50,20){$\sc4$}\put(50,0){$\sc0$}\put(50,-20){$\sc\ol1$} \put
(50,-40){$\sc\ol4$}\put(65,40){$\sc5$} \put(65,20){$\sc2$}\put(65,0){$
\sc\ol2$}\put(65,-20){$\sc\ol3$} \put(65,-40){$\sc\ol6$}\put(80,40){$
\sc2$} \put(80,20){$\sc\ol1$}\put(80,0){$\sc\ol5$}\put(80,-20){$\sc\ol6
$} \put(80,-40){$\sc\ol9$}\put(95,40){$\sc0$} \put(95,20){$\sc\ol3$}\put
(95,0){$\sc\ol7$}\put(95,-20){$\sc\ol8$} \put(95,-40){$\ssc\overline{11}
$}\put(28,-55){$\sc0$}\put(43,-55){$\sc1$} \put(58,-55){$\sc1$}\put
(73,-55){$\sc2$}\put(88,-55){$\sc1$} \put(27,7){\line(3,1){36}}\put
(22,-11){\line(0,1){8}} \put(26,-18){$\sc.$}\put(29,-18){$\sc.$}\put
(32,-18){$\sc.$} \put(27,-33){\line(3,2){21}}\put(40,-13){$\sc.$}\put
(43,-8){$\sc.$} \put(46,-3){$\sc.$}\put(55,-33){$\sc.$}\put(58,-28){$
\sc.$}\put(61,-23) {$\sc.$}\put(52,-10){\line(0,1){7}}\put(70,-11){$\sc.
$}\put(72,-5){$\sc.$} \put(74,1){$\sc.$}\put(76,7){$\sc.$}\put(78,13){$
\sc.$} \put(72,27){\line(3,2){21}}\put(85,27){$\sc.$}\put(88,32){$\sc.
$} \put(91,37){$\sc.$}\put(104,0){\oval(10,76)[r]} \put(115,0){=$\ \ $
\put(0,30){$\sc2$}\put(0,10){$\sc3$}\put(0,-10){$\sc0$}\put(0,-30){$
\sc2$} \put(15,0){\oval(10,76)[l]}\put(20,40){$\sc10$} \put(20,20){$
\sc7$}\put(20,0){$\sc3$}\put(20,-20){$\sc2$}\put(20,-40) {$\sc\ol1$}\put
(35,40){$\sc9$}\put(35,20){$\sc6$}\put(35,0){$\sc2$} \put(35,-20){$\sc1
$}\put(35,-40){$\sc\ol2$}\put(50,40){$\sc7$} \put(50,20){$\sc4$}\put
(50,0){$\sc0$}\put(50,-20){$\sc\ol1$} \put(50,-40){$\sc\ol4$}\put
(65,40){$\sc5$} \put(65,20){$\sc2$}\put(65,0){$\sc\ol2$}\put(65,-20){$
\sc\ol3$} \put(65,-40){$\sc\ol6$}\put(80,40){$\sc2$} \put(80,20){$\sc
\ol1$}\put(80,0){$\sc\ol5$}\put(80,-20){$\sc\ol6$} \put(80,-40){$\sc
\ol9$}\put(95,40){$\put(2,2){\circle{10}}\sc0$} \put(95,20){$\sc\ol3
$}\put(95,0){$\sc\ol7$}\put(95,-20){$\sc\ol8$} \put(95,-40){$\ssc
\overline{11}$}\put(28,-55){$\sc0$}\put(43,-55){$\sc1$} \put(58,-55){$
\sc1$}\put(73,-55){$\sc2$}\put(88,-55){$\sc1$} \put(32,-18){\vector
(-1,0){6}}\put(48,-19){\vector(-3,-2){21}} \put(46,-3){\vector
(-3,-4){8}}\put(52,-3){\vector(0,-1){7}} \put(104,0){\oval
(10,76)[r]}\put(114,0){.}}\par\ni $SW^e_\L(1)=SW_\L(1)=\{(3,3),$ $(4,1),
$ $(4,2),$ $(4,3),$ $(5,1)\}$. However, $SW^e_\L(2)=\{(1,6),$ $(2,4),$ $
(2,5),$ $(3,1),$ $(4,1),$ $(4,4),$ $(5,3)\}$, while $SW_\L(2)=\{(1,6)\}
$. In this case $\g_1,\g_2$ are $n$-related.\qed\par \ni{\sl Example
4.5.} $\L=[1\ 2\ 1\ 1;\ol1;1\ 0\ 0\ 0\ 2]$.\par \ \ \ \ \ \ \ \ \ $A(\L)
$ =$\ \ $\put(0,30){$\sc1$}\put(0,10){$\sc2$} \put(0,-10){$\sc1$}\put
(0,-30){$\sc1$}\put(15,0){\oval(10,76)[l]} \put(20,40){$\sc8$}\put
(20,20){$\sc6$}\put(20,0){$\sc3$}\put(20,-20){$\sc1$} \put(20,-40){$\sc
\ol1$}\put(35,40){$\sc6$}\put(35,20){$\sc4$}\put(35,0) {$\sc1$}\put
(35,-20){$\sc\ol1$}\put(35,-40){$\sc\ol3$}\put(50,40) {$\sc5$}\put
(50,20){$\sc3$}\put(50,0){$\sc0$}\put(50,-20){$\sc\ol2$} \put(50,-40){$
\sc\ol4$}\put(65,40){$\sc4$} \put(65,20){$\sc2$}\put(65,0){$\sc\ol1
$}\put(65,-20){$\sc\ol3$} \put(65,-40){$\sc\ol5$}\put(80,40){$\sc3$}
\put(80,20){$\sc1$}\put(80,0){$\sc\ol2$}\put(80,-20){$\sc\ol4$} \put
(80,-40){$\sc\ol6$}\put(95,40){$\sc0$} \put(95,20){$\sc\ol2$}\put
(95,0){$\sc\ol5$}\put(95,-20){$\sc\ol7$} \put(95,-40){$\sc\ol9$}\put
(28,-55){$\sc1$}\put(43,-55){$\sc0$} \put(58,-55){$\sc0$}\put(73,-55){$
\sc0$}\put(88,-55){$\sc2$} \put(27,-33){\line(3,4){8}}\put(25,-13){$\sc.
$}\put(28,-8){$\sc.$} \put(31,-3){$\sc.$}\put(41,2){$\sc.$}\put(44,2){$
\sc.$}\put(47,2){$\sc.$} \put(42,-13){\line(3,4){8}}\put(57,7){\line
(3,2){21}} \put(56,2){$\sc.$}\put(59,2){$\sc.$}\put(62,2){$\sc.$} \put
(71,2){$\sc.$}\put(74,2){$\sc.$}\put(77,2){$\sc.$} \put(87,27){\line
(3,4){8}} \put(85,9){$\sc.$}\put(87,15){$\sc.$}\put(89,21){$\sc.$}\put
(91,27){$\sc.$} \put(93,33){$\sc.$}\put(102,0){\oval(10,76)[r]}\put
(115,0){=$\ \ $ \put(0,30){$\sc1$}\put(0,10){$\sc2$}\put(0,-10){$\sc1
$}\put(0,-30){$\sc1$} \put(15,0){\oval(10,76)[l]}\put(20,40){$\sc8$}
\put(20,20){$\sc6$}\put(20,0){$\sc3$}\put(20,-20){$\sc1$}\put(20,-40) {$
\sc\ol1$}\put(35,40){$\sc6$}\put(35,20){$\sc4$}\put(35,0){$\sc1$} \put
(35,-20){$\sc\ol1$}\put(35,-40){$\sc\ol3$}\put(50,40){$\sc5$} \put
(50,20){$\sc3$}\put(50,0){$\sc0$}\put(50,-20){$\sc\ol2$} \put(50,-40){$
\sc\ol4$}\put(65,40){$\sc4$} \put(65,20){$\sc2$}\put(65,0){$\sc\ol1
$}\put(65,-20){$\sc\ol3$} \put(65,-40){$\sc\ol5$}\put(80,40){$\sc3$}
\put(80,20){$\sc1$}\put(80,0){$\sc\ol2$}\put(80,-20){$\sc\ol4$} \put
(80,-40){$\sc\ol6$}\put(95,40){$\put(2,2){\circle{10}}\sc0$} \put
(95,20){$\sc\ol2$}\put(95,0){$\sc\ol5$}\put(95,-20){$\sc\ol7$} \put
(95,-40){$\sc\ol9$}\put(28,-55){$\sc1$}\put(43,-55){$\sc0$} \put
(58,-55){$\sc0$}\put(73,-55){$\sc0$}\put(88,-55){$\sc2$} \put
(31,-3){\vector(-3,-4){8}}\put(31,-23){\vector(-3,-4){8}} \put
(47,2){\vector(-1,0){6}}\put(46,-3){\vector(-3,-4){8}} \put(102,0){\oval
(10,76)[r]}\put(112,0){.}}\par\ni $SW^e_\L(1)=SW_\L(1)=\{(3,2),$ $(3,3),
$ $(4,1),$ $(4,2),$ $(5,1)\}$. However, $SW^e_\L(2)=\{(1,6),$ $(2,5),$ $
(3,2),$ $(3,3),$ $(3,4),$ $(3,5),$ $(4,1),$ $(4,2),$ $(5,1)\}$ while $
SW_\L(2)=\{(1,6)\}$. Because $S^e_\L(2),W^e_\L(2)$ meet at the zero
corresponding to $\g_1,SW^e_\L(1)\subset SW^e_\L(2)$, and $\g_1,\g_2$
are $q$-related.\qed\par\ni Although we shall not give details here, it
is easy to check that for doubly atypical $\L=[4020;\ol2;00120],\g_1,\g_
2$ are $n$-related and $W^e_\L(2)$ starts off above $S^e_\L(2)$ but
crosses it at a point which is above the zero corresponding to $\g_1$.
Also, $\L=[1220;\ol2;00110]$ is doubly atypical with $\g_1,\g_2$ being $
q$-related; in this case, $W^e_\L(2)$ starts off above $S^e_\L(2)$ but
meets it at the position of the zero corresponding to $\g_1$; the
position of this zero is therefore not an element of $SW_\L(2)$. These
examples illustrate some properties of chains we are now going to
state.\par\ni {\sl Lemma 4.6.} Let $s<t$. Then $\g_s,\g_t$ are $q
$-related $\Lra$ $W^e_\L(t),S^e_\L(t)$ both contain $(b_s,c_s)$.\vs
{-2pt}\par\ni{\sl Proof.} Since $A(\L)_{b_{\sc s},c_{_{\sc s}}}=A(\L)_
{b_{\sc t},c_{\sc t}}=0$, \vs{-4pt}we have $x_{st}=A(\L)_{b_{\sc t},c_{_
{\sc s}}}=c_t-c_s+\sum^{^{\sc c_{_{\sc t}}-1}}_{_{\sc i=c_{_{\sc s}}}}a_
i$. Also $h_{st}=b_s-b_t+c_t-c_s+1$. \vs{-2pt}Therefore, $\g_s,\g_t$
are $q$-related $\Lra$ $x_{st}=h_{st}-1$ $\Lra$ $b_s=b_t+\sum^{^{\sc c_
{_{\sc t}}-1}}_{_{\sc\,i=c_{_{\sc s}}}}a_i$ $\Lra$ $(b_s,c_s)\in W^e_\L(
t)$. Similarly, if $\g_s,\g_t$ are $q$-related, $(b_s,c_s)\in S^e_\L(t)
$.\qed\par\ni{\sl Lemma 4.7.} If $s<t,\g_s,\g_t$ are $n$-related, then $
W^e_\L(t)$ meets column $c_s$ below $(b_s,c_s)$, and $S^e_\L(t)$ meets
row $b_s$ to the left of $(b_s,c_s)$.\vs{-2pt}\par\ni{\sl Proof.} Since
$\g_s,\g_t$ are $n$-related, $x_{st}>h_{st}-1$, hence\vs{-2pt} $b_s<b_t+
\SUM^{^{\sc c_{_{\sc t}}-1}}_{{\sc\,i=c_{_{\sc s}}}}a_i$. Hence $W^e_\L(
t)$ meets column $c_s$ below $(b_s,c_s)$. The proof of the other result
is analogous.\qed\par\ni {\sl Lemma 4.8.}{\def\SS{{\ssc\!}} If $s\!<\!t,
\g_s,\g_t$ are $q$- or $n$-related, then $SW\!_\L\SS(\SS t\SS)$ does
not extend as far to the left or downwards as $(b_s,\!c_s)$, i.e., $W\!_
\L\SS(\SS t\SS)$ ends to the right of column $c_s$, $S{\sc\!}_\L\SS(\SS
t\SS)$ ends above row $b_s$.} \par\ni{\sl Proof.} If $\g_s,\g_t$ are $q
$-related, by Lemma 4.6 $W^e_\L(t),S^e_\L(t)$ meet at $(b_s,c_s)$.
Otherwise, by Lemma 4.7 $W^e_\L$($t$),$S^e_\L$($t$) must cross to the
right of column $c_s$ and above row $b_s$. In either case, by
Definition 4.2 $SW\!_\L{\sc\!}$(${\ssc\!}t{\ssc\!}$) does not extend as
far to the left or downwards as ($b_s$,$c_s$). \qed\vs{-3pt}\par\ni {\sl
Lemma 4.9.} Let $1\zl t\zl r$ and $b_t\zl d\zl m+1,\ 1\zl e\zl c_t$ and
$i=A(\L)_{de}\in\Z$. If $(d,e)\zi W^e_\L(t)$ then $(d,e+i)\zi S^e_\L(t)
$. If $(d,e)\zi S^e_\L(t)$ then $(d+i,e)\zi W^e_\L(t)$. \vs{-5pt}\par
\ni{\sl Proof.} If $(d,e)\zi W^e_\L(t)$, then $d=b_t+\sum^{^{\sc c_{_{
\sc t}}-1}}_{_{\sc\,j=e}}a_j$. \vs{-3pt}We have $0=A(\L)_{b_{\sc t},c_{_
{\sc t}}}=A(\L)_{de}+\sum^{^{\sc\ol{m-b_{\sc t}+1}}}_{_{\sc j=\ol{m-d+2}
}}a_j-\sum^{^{\sc c_{_{\sc t}}-1}}_{j=e}a_j+d-b_t-c_t+e.$ So $e+i=e+A(
\L)_{de}=c_t-\sum^{^{\sc\ol{m-b_{\sc t}+1}}}_{_{\sc j=\ol{m-d+2}}}a_j$,
\vs{-3pt}by (4.2b), $(d,e+i)\zi S^e_\L(t)$. Similarly, if $(d,e)\zi S^e_
\L(t)$ then $(d+i,e)\zi W^e_\L(t)$. \qed\par\ni {\sl Lemma 4.10.} In
Lemma 4.9, if $(d,e)\zi W_\L(t)$, then $i\ge0$ and $(d,e+i)\zi S_\L(t)
$. If $(d,e)\zi S_\L(t)$, then $i=-j\le0$ and $(d-j,e)\zi W_\L(t)$.\par
\ni{\sl Proof.} Let's say ($d$,$e${}$)\zi W_\L$($t$). By Lemma 4.9 ($d
$,$e$+$i${}$)\zi S^e_\L$($t$). Since ($d$,$e${}$)\zi W_\L$($t$), by
Definition 4.2, $S^e_\L$($t$),$W^e_\L$($t$) must have not crossed each
other above or to the right of ($d$,$e$) and so ($d$,$e$+$i$) must be
to the right of ($d$,$e$). Hence $i${}$\zg$0 and ($d$,$e$+$i$)$\zi${}$S
{\sc\!}_\L$($t$). \qed\par\ni In Example 4.3, $(2,2)\zi W_\L(2),$ $A(\L)
_{22}=3$; one has $(5,2)\zi S_\L(2)$. Also, $(5,1)\zi S_\L(1),$ $A(\L)_
{51}=-1$; one has $(4,1)\zi W_\L(1)$. Finally, $(3,3)\zi W_\L(1)\cap S_
\L(1),\ A(\L)_{33}=0.$ Obviously $W_\L(1)$ is zero rows above $(3,3)$
and $S_\L(1)$ is zero rows to the right of $(3,3)$.\vs{-2pt}\par\ni{\sl
Lemma 4.11.} (i) Let $(d,e)\in W_\L(t),0\zl k\zl m+1-d$. If $A(\L)_{d+k,
e}=j,$ $j+k\zg0$, then $S_\L(t)$ meets the $(d+k)$-th row $(j+k)$
columns to the right of $(d+k,e)$.\nl\hs{3ex}(ii) Let $(d,e)\zi S_\L(t),
0\zl k\zl e-1$. If $A(\L)_{d,e-k}=-j,$ $j+k\zg0$, then $W_\L(t)$ meets
the $(e-k)$-th column $(j+k)$ rows above $(d,e-k)$.\vs{-2pt}\par\ni {\sl
Proof.} Let's prove (i) as (ii) is similar. Let $0\zl k'\zl k,\,A(\L)_
{d+k',e}=j'$, then $j'\zg j+k-k'$ and so $j'+k'\zg j+k\zg0$. The
condition of the lemma is satisfied for $k'$. Let $A(\L)_{de}=i$;\vs
{-3pt} by Lemma 4.10,\vs{-3pt} $i\zg0,\,(d,e+i)\zi S_\L(t)$. By (4.2b),
$(d+k',e+i-\sum_{_{\sc\ell=\ol{m-d+1}}}^{^{\sc\ol{m-d-k'+2}}}a_\ell)\zi
S^e_\L(t).$ \vs{-4pt}Now $A(\L)_{d+k',e}$=$A(\L)_{d,e}${}$-${}$k'${}$-
${}$\sum_{_{\sc\ell=\ol{m-d+1}}}^{^{\sc\ol{m-d-k'+2}}}a_\ell,$ so\vs
{-2pt} $j'$=$i${}$-${}$k'${}$-${}$\sum_{_{\sc\ell=\ol{m-d+1}}}^{^{\sc
\ol{m-d-k'+2}}}a_\ell$ and $(d$+$k',e$+$j'$+$k'$)$\zi${}$S^e_\L(t)$.
Hence, since $j'+k'\zg0$ for all $k',(d+k',e+j'+k')$ is to the right of
column $e$, whereas $W^e_\L(t)$ is to the left of column $e$ after
having passed at $(d,e)$. Thus $(d+k',e+j'+k')\zi S_\L(t).$ So $(d+k,e+
j+k)\zi S_\L(t)$ and the result follows.\qed\par\ni {\sl Definition
4.12.} Suppose $W_\L(t),S_\L(t)$ end at $(d_t,e_t),(d'_t,e'_t)$,
respectively, $d'_t\zg d_t,e'_t\zg e_t$. Define $D(t)$ to be the region
of $D$ within or on the boundary consisting of $SW_\L(t)$, the vertical
line joining $(d'_t,e_t)$ to $(d_t,e_t)$ and the horizontal line
joining $(d'_t,e_t)$ to $(d'_t,e'_t)$.\qed\vs{4pt}\nl We will prove
that every element of $D(t)$ lies in some $SW_\L(s)$ such that $s${}$
\zl${}$t,$ $\g_s,\g_t$ are $c$-related.\vs{-2pt}\par\ni {\sl Lemma
4.13.} (i) If $A(\L)_{bc}>0$ and $A(\L)_{b-a_c,c+1}\zl 0$ then $A(\L)_
{bc}=1$ and $A(\L)_{b-a_c,c+1}=0.$\nl\hs{3ex}(ii) If $A(\L)_{bc}<0$ and
$A(\L)_{b-1,c+a_{\ol{m-b+2}}}\zg0$ then $A(\L)_{bc}=-1$ and $A(\L)_{b-1,
c+a_{\ol{m-b+2}}}=0.$\vs{-2pt} \par\ni{\sl Proof.} The proof of (ii) is
analogous to that of (i). For (i), using (3.1a),\vs{-10pt} $$A(\L)_{b-a_
c,c+1}=A(\L)_{bc}-(a_c+1)+\rb{3pt}{\mbox{$ _{^{\ \ \ \ \sum}_{\sc\ell
=\ol{m-b+2}}}^{\ol{m- b+1+ a_c}} $}}a_\ell+a_c=A(\L)_{bc}-1+\rb{3pt}{
\mbox{$ _{^{\ \ \ \ \sum}_{\ell=\ol{m- b+2}}}^{\ol{m-b+1+ a_c}} $}}a_
\ell\zl0.$$\vs{-14pt}\nl The only solution is $a_{_{\sc\ol{m-b+2}}}$=$
\sc\cdots$=$a_{_{\sc\ol{m-b+1+a_c}}}$=$0,A(\L)_{bc}$=$1,A(\L)_{b-a_c,c+
1}$=0.\qed\vs{-2pt}\par\ni {\sl Lemma 4.14.} Let $(d,e)\in D(t)$. \vs
{-1pt}\nl\hs{3ex}(i) If $A(\L)_{de}\zg0$, then there exists $s,\,1\zl s
\zl t$ such that $(b_s,c_s)\zi D(t),\,(d,e)\zi W_\L(s)$; \vs{-1pt}\nl\hs
{3ex}(ii) If $A(\L)_{de}\zl0$, then there exists $s,\,1\zl s\zl t$ such
that $(b_s,c_s)\zi D(t),\,(d,e)\zi S_\L(s)$; \vs{-1pt}\nl\hs{3ex}(iii)
In both of these cases, $\g_s,\g_t$ are $c$-related. \vs{-2pt}\par\ni
{\sl Proof.} (i) Suppose $A(\L)_{de}\zg0$ and let $k$ be such that $(d-
k,e)\zi W_\L(t)$, where clearly $0\zl k\zl d-1$. If $k=0$, then $(d,e)
\zi W_\L(t)$.\vs{-3pt} Suppose $k>0$. Define the ordered set $W$ of
elements of $D$ by $W$=$\{(d$,$e),(d${}$-${}$a_e$,$e$+1),$\cdots$,($d
${}$-${}$\sum_{_{\sc i=e}}^{^{\sc c_{_{\sc t}}-1}}a_i$,$c_t$)\}.\vs
{-2pt} Each element of $W$ is $k$ rows below an element of $W_\L(t)$,
so $W$ must meet or cross $S_\L(t)$, by which time the corresponding
entries of $A(\L)$ have become negative. By Lemma 4.13 there must be an
element of $W$ for which the corresponding entry of $A(\L)$ is zero,
this element must lie to the left of $S_\L(t)$, hence in $D(t)$, i.e.,
there is $\g_s,1${}$\zl${}$s${}$\zl${}$t$, such that $(b_s$,$c_s$)$\zi
${}$W${}$\cap${}$D(t)$ and so $(d$,$e$)$\zi${}$W_\L(s)$. The proof of
(ii) is analogous to (i). For (iii), in both cases, $(b_s$,$c_s$)$\zi
${}$D(t)$; suppose $(b_s$,$c_s)$ is\vs{-3pt} $k$ rows below $W_\L(t)$,
i.e., $(b_s${}$-${}$k$,$c_s$)$\zi${}$W_\L(t)$, $k${}$\zg$1. We\vs{-8pt}
then have $b_s${}$-${}$k$=$b_t$+$\sum_{_{\sc i=c_{_{\sc s}}}}^{^{\sc c_
{_{\sc t}}-1}}a_i.$ Now, since $A(\L)_{b_{\sc t},c_{_{\sc t}}}$=0, we
have $x_{st}$=$\sum_{_{\sc i=c_{_{\sc s}}}}^{^{\sc c_{_{\sc t}}-1}}a_i
$+$c_t${}$-${}$c_s$=$b_s${}$-${}$b_t$+$c_t${}$-${}$c_s${}$-${}$k$=$h_
{st}${}$-$1$-${}$k${}$<${}$h_{st}${}$-$1. Hence $\g_s,\g_t$ are $c
$-related.\qed\vs{-2pt}\par \ni{\sl Lemma 4.15.} If $1\zl s<t\zl r$ and
$\g_s,\g_{s+1},{\sc\ldots},\g_{t-1}$ are all $c$-related to $\g_t$,
then $D(s)\subset D(t)$.\vs{-2pt}\par\ni{\sl Proof.} Let $b$=$b_t$+$\vs
{-6pt}\sum_{_{\sc i=c_{_{\sc s}}}}^{^{\sc c_{_{\sc t}}-1}}a_i$, by
(4.2a), $(b$,$c_s$)$\in${}$W^e_\L(t)$. Since $h_{st}$=$c_t${}$-${}$c_s
$+$b_s${}$-${}$b_t$+1, $x_{st}$=$c_t${}$-${}$c_s$+ $\sum_{_{\sc i=c_{_{
\sc s}}}}^{^{\sc c_{_{\sc t}}-1}}a_i,$ so $b=b_t$+$x_{st}${}$-${}$c_t
$+$c_s<b_t$+$h_{st}${}$-$1$-${}$c_t$+$c_s=b_s$.\vs{-6pt} i.e., $(b$,$c_
s)$ is above $(b_s$,$c_s)$ and since $A(\L)_{b_{\sc s},c_{_{\sc s}}}
$=0, so the entry $A(\L)_{b,c_{_{\sc t}}}${}$>$0. \vs{-3pt}Similarly,
let $c$=$c_s${}$-${}$\sum^{^{\sc b_{\sc s}-1}}_{_{\sc i=b_{\sc t}}}a_{_
{\sc\ol{m-i+1}}}$, then $(b_s$,$c$)$\in${}$S^e_\L(t)$ is to the right
of $(b_s$,$c_s)$ and the entry $A(\L)_{b_{\sc s},c}${}$<$0. Elements $(
b,c_s)$, $(b_s,c)$ will be in fact in $W_\L(t),S_\L(t)$, respectively,
provided $W^e_\L(t),S^e_\L(t)$ do not cross above and to the right of
these elements. Suppose they do cross; then the entries corresponding
to elements of $W^e_\L(t)$, $S^e_\L(t)$ must change from, respectively,
positive to negative, negative to positive above and to the right of
these elements. By Lemma 4.13, this means that $W^e_\L(t),S^e_\L(t)$
contain a common element $(d,e)$ with $A(\L)_{de}=0$, i.e., there is $
\g_p$, $s<p<t$ such that $(d,e)=(b_p,c_p)$. By Lemma 4.6, $\g_p,\g_t$
must be $q$-related, contradiction with the assumption of the lemma.
This proves that ($b$,$c_s$)$\in${}$W_\L$($t$), ($b_s$,$c$)$\in${}$S_\L(
t)$. It follows that $(b_s,c_s)\zi D(t).$ From this, it follows that $W^
e_\L(s),S^e_\L(s)$ must cross before $W^e_\L(t),S^e_\L(t)$ cross, so $
SW_\L$($s$)$\subset${}$D$($t$).\,Hence\,also\,$D$($s$)$\subset${}$D$($t
$).\qed\par\ni Lemmas 4.14-5 are illustrated by Example 4.3, where $r=2,
\g_1,\g_2$ are $c$-related. Clearly, $D(1)\subset D(2)$ and every
position of $D(2)$ are in $SW_\L(2)$ or $SW_\L(1)$. Now the definition
below tells how to determine the composition factors of $\VBL$
corresponding to the unlinked codes.\par\ni {\sl Definition 4.16.} Let $
\Si^c$ be an unlinked code for $\L$ with non-zero columns $C_{s_1},{\sc
\ldots},C_{s_p}$ where $1\zl s_1<{\sc\ldots}<s_p\zl r$. Define the
subset $D_\Si$ of $D$ to be $\phi$ if $p=0$ or $\cup_{l=1}^p SW_\L(s_l)
$ otherwise. The weight corresponding to the code $\Si^c$ is then
defined by $\Si=\L-\sum_{\b\in\wh D_\Si}\b.$\qed\vs{3pt}\par\ni $D_\Si,
\Si$ are uniquely determined by $\Si^c$. If $\Si^c$=0${\sc\ldots}$0,
then $D_\Si$=$\phi,\Si$=$\L$. Note that not all unions of chains
correspond to codes. Recall that if, for all $u$, $s${}$\zl${}$u${}$<
${}$t,$ $\g_u$ is $c$-related to $\g_t$ then in all codes $\Si^c$ with
non-zero column $t$ its $s$-th column must contain $t$, i.e., $\g_s$
must be wrapped by $\g_t$. Hence if $SW_\L(t)${}$\subset${}$D_\Si$,
then also $SW_\L(u)${}$\subset${}$D_\Si$ for all $u$, $s${}$\zl${}$u${}$
<${}$t$. Now if $s$ is the smallest number such that for all $u$, $s${}$
\zl${}$u${}$<${}$t,$ $\g_u$ is $c$-related to $\g_t$, then Lemmas
4.14-5 show that $D(t)$=$\cup^t_{u=s}SW_\L(u)$. Thus the requirement on
codes $\Si^c$ with non-zero column $t$ that $\g_t$ must wrap all $\g_u
$, $s${}$\zl${}$u${}$<${}$t$ is equivalent to the requirement that $D(t)
${}$\subset${}$D_\Si$. A union of south west chains not satisfying this
condition cannot correspond to a code. If $\Si^c$ is an indecomposable
unlinked code, then each non-zero column of the code contains a common
number, $t$ say. This means that column $t$ is the rightmost non-zero
column of $\Si^c$ and $\g_t$ wraps all the $\g$'s corresponding to the
other non-zero columns, so that $D_\Si$=$D(t)$. Conversely, if $D_\Si$=$
D(t)$ for some $t$, we may reverse the above argument to show $\Si^c$
is indecomposable. We therefore have\par\ni {\bf Theorem 4.17.} $\Si^c$
is an indecomposable unlinked code \ $\Lra$ \ there exists a $t,D_\Si=D(
t)$.\qed\par \ni{\sl Example 4.18.} In the case $sl(6/5)$ let $\L=[
00020;0;0210]$ with $A(\L)$ and $nqc(\L)$ as given in (3.3). The south
west chains of $A(\L)$ are as follows:\par\hs{90pt}\ \ \ \ \ \ \ \ \ \ $
A(\L)$ =$\ \ $\put(0,40){$\sc0$}\put(0,20){$\sc0$}\put(0,0){$\sc0$}\put
(0,-20) {$\sc2$}\put(0,-40){$\sc0$}\put(15,0){\oval(10,105)[l]}\put
(20,50){$\sc7$} \put(20,30){$\sc6$}\put(20,10){$\sc5$}\put(20,-10){$
\sc4$}\put(20,-30) {$\sc1$}\put(20,-50){$\put(2,2){\circle{10}}\sc0
$}\put(40,50){$\sc6$} \put(40,30){$\sc5$}\put(40,10){$\sc4$}\put
(40,-10){$\sc3$}\put(40,-30) {$\sc0$}\put(40,-50){$\sc\ol1$}\put
(60,50){$\sc3$}\put(60,30){$\sc2$} \put(60,10){$\sc1$}\put(60,-10){$
\put(2,2){\circle{10}}\sc0$}\put(60,-30) {$\sc\ol3$}\put(60,-50){$\sc
\ol4$}\put(80,50){$\sc1$}\put(80,30){$\sc0$} \put(80,10){$\sc\ol1$}\put
(80,-10){$\sc\ol2$}\put(80,-30){$\sc\ol5$} \put(80,-50){$\sc\ol6$}\put
(100,50){$\sc0$}\put(100,30){$\sc\ol1$} \put(100,10){$\sc\ol2$}\put
(100,-10){$\sc\ol3$}\put(100,-30){$\sc\ol6$} \put(100,-50){$\sc\ol7
$}\put(110,0){\oval(10,105)[r]}\put(30,-65){$\sc0$} \put(50,-65){$\sc2
$}\put(70,-65){$\sc1$}\put(90,-65){$\sc0$} \put(38,-8){\vector
(-1,0){8}}\put(38,-28){\vector(-1,0){8}} \put(42,-32){\vector
(0,-1){8}}\put(58,28){\vector(-1,-2){15}} \put(62,-32){\vector
(0,-1){8}}\put(78,48){\vector(-1,-1){11}} \put(78,28){\vector
(-1,-1){11}}\put(82,28){\vector(0,-1){8}} \put(82,8){\vector
(0,-1){8}}\put(98,52){\vector(-1,0){8}} \put(98,-12){\vector
(-2,-1){30}}\put(102,48){\vector(0,-1){8}} \put(102,28){\vector
(0,-1){8}}\put(102,8){\vector(0,-1){8}}\hs{120pt}. \hfill(4.3)\par\ni
All codes were given in (3.4). One decomposable unlinked code is $\Si^c
$=${\ssc\!}\,^{1\,0\,}{}^{3\,4\,0}_4$, in which $\g_4$ wraps $\g_3$. \vs
{-3pt}The non-zero columns are the 1st, 3rd, 4th and $D_\Si=SW_\L(1)
\cup SW_\L(3)\cup SW_\L(4)=D(1)\cup D(4)$, $\Si=\L-\b_{61}-\b_{43}-(\b_
{33}+\b_{24}+\b_{34}+\b_{44})$. \vs{-7pt}One indecomposable unlinked
code is \rb{-8pt}{$_{\,}^{^{^{\sc1\,2}_{\sc2\,5}}_{\sc5}\,^{^{\sc3\,4
\,5}_{\sc4\,5}}_{\sc5}}$}, in which $\g_2$ wraps $\g_1,$ $\g_4$ wraps $
\g_3$ and $\g_5$ wraps all $\g_1,{\sc\ldots},\g_4$ and $D_\Si=\cup_{s=1}
 ^5SW_\L(s)=D(5)$ and $\Si=\L-\b_{61}-(\b_{51}+\b_{52}+\b_{62})-\b_{43}-(
\b_{33}+\b_{24}+\b_{34}+\b_{44})-(\b_{41}+\b_{42}+\b_{23}+\b_{14}+\b_
{15}+\b_{25}+\b_{35}+\b_{45}+\b_{53}+\b_{63})$. On the other hand, $SW_
\L(2)$ and $SW_\L(1)\cup SW_\L(3)\cup SW_\L(5)$ do not correspond to
codes since they violate the necessary conditions that $\g_2$ must wrap
$\g_1$ and $\g_5$ must wrap all of $\g_1,{\sc\ldots},\g_4$.\qed\par\ni
It was part of the conjecture$^{\,2}$ that for all unlinked codes $\Si^
c$ of $\L$, $\Si$ defined in Definition 4.16 is the highest weight of a
composition factor $V(\Si)$ of $\VBL$ (note that in {\sl Ref.} 2, $\Si$
was defined by means of boundary strip removals; for unlinked codes it
is easy to see that this is equivalent to the definition here). Note
that this implies that $\Si$ is a dominant weight, which was in fact
proved using the corresponding strip removals in the Young diagram. For
$\Si^c=0{\sc\ldots}0,\Si=\L$, we see that this code corresponds to the
top composition factor $V(\L)$. We shall prove that, for any unlinked
code $\Si^c$ for $\L$, there exists a primitive vector $v_{\ssc\Si}$
and correspondingly a composition factor $V(\Si)$ of $\VBL$. To make
connection with codes explicit, we use notation $v_{\ssc\Si^c},$ $U(\Si^
c),$ $V(\Si^c)$ to denote, respectively, $v_{\ssc\Si}$, ${\bf U}(G)v_{
\ssc\Si}$, $V(\Si)$. Thus, if $\Si^c=0{\sc\ldots}0$, then $v_{(0
\ldots0)}\equiv v_{\ssc\L},U(0{\sc\ldots}0)\equiv\VBL$ and $V(0{\sc
\ldots}0)\equiv V(\L)$. Finally, it is worth mentioning that the linked
codes for $\L$, if any occur, appear to correspond to south west chains
in $A(\Pi^+)$, where $\Pi^+=\Pi+2\rho_1$ and $\Pi$ is the lowest $G_{
\ol0}$-highest weight of the simple module $V(\Si)$ corresponding to
the code; more details of this are given in {\sl Ref.} 2.\es {\bf V.
MORE NOTATION AND PRELIMINARY RESULTS}\vs{-2pt}\par\ni \vs{-1pt}Define
a total order on $\D$: $\a_{ij}<\a_{k\ell}\Lra j-i<\ell-k$ or $j-i=\ell-
k$, $i>k.$ It implies that $\b_{bc}<\b_{de}\Lra c-b<d-e$ or $c-b=d-e$, $
b>d.$ This total order on\vs{-1pt} $\D_1$ corresponds to the sequence
of positions signified by $1,2,{\sc\cdots}$ in the following $(m+1)
\times(n+1)$ matrix, where $\b_{bc}$ is the root associated with the $(
b,c)$-th entry:\nl\hs{180pt} $\left(\mbox{\rb{4pt}{$\matrix{\sc\cdot\
\,\cdot\ \,\cdot\ \cdots \hfill\vs{-7pt}\cr \sc6\ \cdot\ \,\cdot\ \cdots
\hfill\vs{-5pt}\cr \sc3\ 5\ \cdot\ \cdots\hfill\vs{-5pt}\cr \sc1\ 2\ 4\
\cdots\vs{-5pt}\cr}$}}\right)$\hfill(5.1)\vs{2pt}\nl \vs{-1pt}By
Theorem 2.1, choose a \vs{-2pt}basis $B$ of ${\bf U}(G_{-1})$: $B=\{b={
\sc\prod}_{\b\in S}f(\b)|S\subset\D^+_1\},$ where $f(\b)$ is a
negative\vs{-1pt} root vector corresponding to $\b$ and the product ${
\sc\prod}_{\b\in S}f(\b)=f(\b_1){\sc\cdots}f(\b_s)$ is written in the
{\sl proper order}\,: $\b_1<{\sc\cdots}<\b_s$ and $s=|S|$ (the {\sl
depth} of $b$). Define a total order on $B$:\vs{-4pt} $$b>b'=f(\b'_1){
\sc\cdots}f(\b'_{s'})\ \Lra\ s>s'\mbox{ \ or \ }s=s',\b_k>\b'_k,\,\b_i=
\b'_i\ (1\le i\le k-1),$$\vs{-18pt}\nl where $b,b'$ are in proper
order. Recall that an element $v\!\in\!\VBL$ can be uniquely written
as\vs{-4pt} $$v=b_1y_1v_\L+b_2y_2v_\L+{\sc\cdots}={\sc\sum}^t_{i=1}b_iy_
iv_\L,\ b_i\zi B,b_1>b_2>{\sc\cdots},\ 0\ne y_i\zi{\bf U}(G_0^-).\eqno(
5.2)$$\vs{-14pt}\nl Clearly\,$v$=0$\Lra${}$t$=0. If\,$v${}$\ne
$0,\,we\,call\,$b_1y_1v_\L$\,the\,{\sl leading\,term}. A\,term\,$b_iy_
iv_\L$\,is\,called\,a\,{\sl prime\,term} if $y_i${}$\in${}$\C
$.\,Note\,that\,a\,vector\,$v
$\,may\,have\,zero\,or\,more\,than\,one\,prime
terms.\,One\,immediately\,has\par\ni {\sl Lemma 5.1.} Let $v$=$gu$, $u
${}$\in${}$\VBL$, $g${}$\in${}${\bf U}(G^-)$. (i)\,If\,$u
$\,has\,no\,prime\,term then $v$ has\,no\,prime\,term. (ii)\,Let\,$v'$=$
gu',u'${}$\in${}$\VBL$.\,If\,$u,u'
$\,have\,the\,same\,prime\,terms\,then\,$v,v'
$\,have\,the\,same\,prime\,terms.\qed\par\ni {\sl Lemma 5.2.} (i) Let $
v_{\ssc\Si}\in\VBL$ be a $G_{\ol0}$-primitive vector with weight $\Si$.
Then $\L-\Si$ is a sum of distinct positive odd roots, furthermore the
leading term $b_1y_1v_\L$ of $v_{\ssc\Si}$ must be a prime term.\nl\hs
{3ex}(ii) Suppose $v'_{\ssc\Si}=\sum^{t'}_{i=1}(b'_iy'_i)v_\L$ is
another $G_{\ol0}$-primitive vector with weight $\Si$. If all prime
terms of $v_{\ssc\Si}$ are the same as those of $v'_{\ssc\Si}$, then $v_
{\ssc\Si}=v'_{\ssc\Si}$.\par\ni{\sl Proof.} (i) Let $v_{\ssc\Si}$ be as
in (5.2). If $y_1\notin\C$, then there exists $e_k\zi G_0^+,$ $(e_ky_1)
v_\L\ne0.$ We have\vs{5pt}\nl\cl{ $e_kv_{\ssc\Si}=b_1(e_ky_1)v_\L+[e_k,
b_1]y_1v_\L+b_2(e_ky_2)v_\L+[e_k,b_2]y_2v_\L+\cdots,$}\vs{6pt}\nl and $[
e_k,b]={\sc\sum}_{p=1}^s f(\b_1){\sc\cdots}[e_k,f(\b_p)]{\sc\cdots}f(\b_
s)$ if $b=f(\b_1){\sc\cdots}f(\b_s)$ and (2.4) gives $[e_k,f(\b_p)]$ =0
or $\pm f(\b_q)$ with $\b_q<\b_p$. The leading term of $b_i(e_ky_i)v_\L+
[e_k,b_i]y_iv_\L$ is $b_i(e_ky_i)v_\L$ if $e_ky_i${}$\ne$0. It follows
that $b_1(e_ky_1)v_\L$ is the leading term of $e_kv_{\ssc\Si}$, i.e., $
e_kv_{\ssc\Si}\zn0$, contradicting that $v_{\ssc\Si}$ is $G_{\ol0}
$-primitive. So, $y_1\zi\C$ and $\L-\Si$ is the weight of $b_1$, a sum
of distinct positive odd roots.\vs{-4pt}\par(ii) Let $v=v_{\ssc\Si}-v'_
{\ssc\Si}$. If $v\zn0$ (then it must be $G_{\ol0}$-primitive), since
its prime terms are all cancelled, $v$ has no prime term, therefore by
(i), it is not $G_{\ol0}$-primitive, a contradiction. \qed\par\ni{\sl
Lemma 5.3.} Let $\Si$ be a (weakly) primitive weight of $\VBL$. Then $
\la\L+\rho|\L+\rho\ra=\la\Si+\rho|\Si+\rho\ra.$\par\ni {\sl Proof.}
Using the Casimir operator $\O=2\u^{-1}(\rho)+{\sc\sum}_i u^iu_i+2{\sc
\sum}_{\a\in\D^+}f(\a)e(\a)$, where $\{u^i\}$ is a basis of $H,\{u_i\}$
is its deal basis, $\u$ is the isomorphism: $H$* $\rar H$ derived from $
\la\cdot|\cdot\ra$, {\sl cf.} that of Lie algebras,$^7$ we see that $\O|
_{\VBL}=\la\L+2\rho|\L\ra I|_{\VBL}.$ Hence, since $\Si$ is the weight
of a weakly primitive vector, we have $\la\L+2\rho|\L\ra=\la\Si+2\rho|
\Si\ra$.\qed\vs{-2pt}\par\ni Define the following zero weight elements
of ${\bf U}(G_{\ol{0}})$:\vs{-10pt} $$\o_i\!=\!\cases{m+1+ i+ h_{\ol mi}
, &$i\zi I_1$,\cr n+1- i+ h_{in},   &$i\zi I_2$,\cr}\ \ \ \ \ \O_i\!=\!
\cases{1-{\sc\sum}_{\ol m\le\ol k\le i}f_{\ol ki}e_{\ol ki}, & $i\zi
I_1$,\cr 1-{\sc\sum}_{i\le k\le n}f_{ik}e_{ik},   &$i\zi I_2$,\cr}\ \ \
\ \ X_i\!=\!\o_i+\O_i.\eqno(5.3)$$\vs{-5pt}\nl {\sl Lemma 5.4.} $[\o_i,
\o_j]=[\o_i,\O_j]=[\O_i,\O_j]=[X_i,X_j]=0$ for all $i,j\in I_1\cup I_2.
$\par\ni{\sl Proof.} $\o_i\in\C\oplus H$ and $\O_i$ with weight 0 imply
the vanishing of the first 2 commutators. The 3rd vanishes trivially if
$i\in I_1,$ $j\in I_2$ or {\it vice versa}. Say, $i,j\in I_1,$ $i<j$.
For each summand of $\O_i$ we have $[f_{ki}e_{ki},\O_j]=-[f_{ki}e_{ki},
f_{kj}e_{kj}+f_{i+1,j}e_{i+1,j}+\cdots]=0$ for $\ol m\le k\le i\le j\le
\ol1,$ where, by (2.4), the omitted terms commute with $f_{ki}e_{ki}$
and 2 non-vanishing terms are cancelled. Thus the 3rd vanishes. The
vanishing of the first three implies the vanishing of the 4th.\qed\vs
{-1pt}\par\ni Now, from the definition of ${\bf U}(G)$ and the $\Z
$-grading of $G$, we can define a projection:\vs{-6pt} $$\v:\,{\bf U}(G)
\rar{\bf U}(G^-\oplus H)\mbox{ derived from }{\bf U}(G)={\bf U}(G^-
\oplus H)\oplus{\bf U}(G)G^+,\eqno(5.4)$$ where ${\bf U}(G)G^+$ is a
left ideal of ${\bf U}(G)$. For $i\in I_1\cup I_2,c\in\C$ and $g\in{\bf
U}(G^-\oplus H)$ define $\chi_{i,c},\chi_i$: ${\bf U}(G^-\oplus H)\rar{
\bf U}(G^-\oplus H)$ by (where $\equiv$ means equal under \ mod\,${\bf
U}(G)G^+$)\vs{-4pt} $$\chi_{i,c}g=\v((c+\O_i)g)=cg+\v(\O_ig)\equiv cg+
\O_i g,\ \ \chi_ig=\v(X_ig)=\v((\o_i+\O_i)g)\equiv X_ig.\eqno(5.5)$$
{\sl Lemma 5.5.} The operators $\chi_{i,c}$ and $\chi_{j,c'}$ commute,
so do $\chi_i$ and $\chi_j$.\par\ni{\sl Proof.} For $g\in{\bf U}(G^-
\oplus H),$ $\chi_{i,c}\chi_{j,c'}g$=$\chi_{i,c}(c'g+\v(\O_jg))$=$c(c'g+
\v(\O_jg))+\v(\O_i(c'g+\v(\O_jg)))$=$cc'g+c\v(\O_jg)+c'\v(\O_ig)+\v(\O_
i\v(\O_jg))$ and $\v(\O_i\v(\O_jg))\equiv\O_i\v(\O_jg)\equiv\O_i\O_jg$
mod$\,{\bf U}(G)G^+$ (the 1st formula follows from (5.4), the 2nd from
the fact that ${\bf U}(G)G^+$ is a left ideal of ${\bf U}(G)\,$). This
gives the 1st part of the lemma by virtue of the commutativity of $\O_i,
\O_j$. Similarly, $\chi_{i}\chi_{j}g\!=\!\chi_i\v(X_jg)=\v(X_i\v(X_jg))
\equiv X_i\v(X_jg)\equiv X_iX_jg$ mod$\,{\bf U}(G)G^+$ and the
commutativity of $X_i,X_j$ implies the 2nd part of the lemma.\qed\par\ni
This Lemma allows us to make the following definitions. For any $J=\{j_
1,j_2,\ldots\}\subseteq I_1\cup I_2,C=(c_{j_1},c_{j_2},\ldots)\in\C^{
\otimes\#J},g\in{\bf U}(G^-\oplus H)$, let\vs{3pt}\nl\hs{50pt} $\chi_{
\ssc J,C}g=\pr{\,j\in J}\chi_{{\ssc j},c_j}g={\sc\cdots}\chi_{{\ssc j_2}
,c_{j_2}}\chi_{{\ssc j_1},c_{j_1}}g,\ \ \ \ \ \ \chi_{\ssc J}g=\pr{\,j
\in J}\chi_{\ssc j}g={\sc\cdots}\chi_{\ssc j_2}\chi_{\ssc j_1}g.$ \hfill
(5.6)\vs{4pt}\nl Now we are in a position to establish some important
results for the successive application of $\chi_{i,c_i}$ to $f_{\ol{r},
s}$ for various $i\in I_1\cup I_2$. From (5.3\&5.5) and (2.4), for $\ol
m\le\ol r\le0\le s\le n$, we have $$\chi_{\ol i,c_{\ol i}}f_{\ol rs}=c_
{\ol i}f_{\ol rs}+f_{\ol{i-1},s}f_{\ol r\ol i},\,\ol r\!\le\!\ol i\!<
\!0\mbox{ \ and \ }\chi_{i,c_i}f_{\ol rs}=c_if_{\ol rs}-f_{\ol r,i-1}f_
{is},\,0\!<\!i\!\le\!s.\eqno(5.7)$$ Further application of the
commutation relations gives:\par\ni\cl{ $\matrix{\chi_{\ol j,c_{\ol j}}
\chi_{\ol i,c_{\ol i}}f_{\ol rs}=c_{\ol j}c_{\ol i}f_{\ol rs}+c_{\ol j}
f_{\ol{i-1},s}f_{\ol r\ol i}+c_{\ol i}f_{\ol{j-1},s}f_{\ol r\ol j}+f_{
\ol{i-1},s}f_{\ol{j-1},\ol i}f_{\ol r\ol j},&\ol r\zl\ol j<\ol i<0,
\hfill\cr\chi_{\ol j,c_{\ol j}}\chi_{i,c_i}f_{\ol rs}=c_{\ol j}c_if_{
\ol rs}-c_{\ol j}f_{\ol r,i-1}f_{is}+c_if_{\ol{j-1},s}f_{\ol r\ol j}-f_
{\ol{j-1},i-1}f_{\ol r\ol j}f_{is},&\ol r\!\zl\!\ol j\!<\!0\!<\!i\!\zl
\!s,\hfill\cr\chi_{j,c_j}\chi_{i,c_i}f_{\ol rs}=c_jc_if_{\ol rs}-c_jf_{
\ol r,i-1}f_{is}-c_if_{\ol r,j-1}f_{js}+f_{\ol r,i-1}f_{i,j-1}f_{js},&0<
j<i\zl s.\hfill\cr}$ } \par\ni The pattern of terms is becoming clear.
The following result may be proved inductively: \par\ni{\sl Lemma 5.6.}
Let $J=\ol P\cup Q$ with $\ol P\subseteq\{\ol r,\ldots,\ol1\},\ Q
\subseteq\{1,\ldots,s\},\ 1\le r\le m,\,1\le s\le n$. Then, with the
definition (5.6),\vs{-2pt}\nl\hs{30pt} $\matrix{\chi_{\ssc J,C}f_{\ol
rs}=&\rb{4pt}{\mbox{$ _{^{\sum}_{k=0}}^{\#\ol P}\ _{^{\sum}_{\ell
=0}}^{\#Q} $}}\ (-1)^\ell\ \ \rb{-3pt}{\mbox{$ ^{_{\sc\ \ \ \sum
}}_{^{\sc\ol K=\{\ol{j_k},\ldots,\ol{j_1}\}\subseteq\ol P} _{\sc\ \ \ol
{j_k}<\cdots<\ol{j_1}}}\ \ ^{_{\sc\ \ \ \sum}}_{^{\sc L=\{i_1,\ldots
,i_\ell\}\subseteq Q} _{\sc\ \ i_1<\cdots<i_\ell}}$}}\ \pr{\sc\ol j\in
\ol P\bs\ol K}c_{\ol j}\ \pr{\sc i\in Q\bs L}c_i\,\cdot\hfill\vs{3pt}
\cr&\hfill f_{\ol{j_0}i_0}\cdot f_{\ol{j_2-1},\ol{j_1}}f_{\ol{j_3-1},
\ol{j_2}}{\sc\cdots}f_{\ol r\ol{j_k}}\cdot f_{i_1,i_2-1}f_{i_2,i_3-1}{
\sc\cdots}f_{i_\ell s}\cr}$\hfill(5.8)\vs{-2pt}\par\ni where $\ol{j_0}
$=$\ol r$ or $\ol{j_1\!-\!1}$ if $k$=$0$ or not; $i_0$=$s$ or $i_1\!-
\!1$ if $\ell$=$0$ or not. Similarly, successive application of $\chi_i
$ to $f_{\ol{r}s}$ gives results exactly analogous to those of (5.7-8)
with $c_i$ replaced by $\o_i$.\qed\par\ni In (5.8) negative root
vectors $f_{ij}$ correspond to $\a_{ij}${}$\zi${}$\D^+$ and the
products of root vectors have been ordered in such a way that the
leftmost factor $f_{\ol{j_0}i_0}$ is a odd vector, while the remaining
factors $f_{ij}$ are even. Moreover, in every summand the elements $c_{
\ol j},c_i$ or $\o_{\ol j},\o_i$, which lie in $\C${}$\oplus${}$H$,
precede an element (a product of $f_{ij}$) of ${\bf U}(G^-)$ which in
every case have weight $-\a_{{\ol r}s}$. From this follows the crucial
relationship linking $\chi_{\ssc J}$ and $\chi_{\ssc J,C}$. For any
weight $\l$, define $$c_i(\l)=\cases{\sum^i_{_{\sc k=\ol m}}\l(h_k)+m+i&
for $i\in I_1$, \vs{2pt}\cr \sum^n_{k=i}\l(h_k)+n-i&for $i\in I_2$.\cr}
\eqno(5.9)$$ With this notation and Lemma 5.5 in the special case for
which $r=m$ and $s=n$, we have\par\ni{\sl Corollary 5.7.} Let $v_\l$
have weight $\l$. Then $\chi_{\ssc J}f_{\ol mn}v_\l=\chi_{\ssc J,C}f_{
\ol mn}v_\l$ with $c_i=c_i(\l),\,i\zi J.$\par\ni {\sl Proof.} It
follows from the above remarks about the order and nature of the
factors in (5.8) that for each $i\in J$ the factor $\o_i$, defined by
(5.3), gives rise to a factor $c_i$ in (5.8) with $$c_i=\cases{m+1+i+\l(
h_{\ol mi})-\a_{\ol mn}(h_{\ol mi}) &for $i\in\ol P\subseteq I_2$,\cr n+
1-i+\l(h_{in})-\a_{\ol mn}(h_{in})& for $i\in Q\subseteq I_1$.\cr}$$
But\vs{-3pt} $\a_{\ol mn}(h_j)=1$ if $j=\ol m$ or $n$; $=0$ otherwise
and $h_{\ol mi}=\sum^i_{_{\sc k=\ol m}}h_k,\,i\zi I_1$ and\vs{-3pt} $h_
{in}=\sum^n_{k=i}h_k,\,i\zi I_2$. It follows that $\a_{\ol mn}(h_{\ol
mi})=\a_{\ol mn}(h_{in})=1$ so that $c_i=c_i(\l)$ as required. \qed\par
\ni It\vs{-3pt} is also worth observing that the explicit expansion
(5.8) implies: \par\ni{\sl Corollary 5.8.} With notation of Lemma 5.6, $
\chi_{\ssc J,C}f_{\ol rs}=\chi_{\ssc J,C}^{(r/s)}f_{\ol rs}$.\qed\par\ni
We\vs{-3pt} shall need commutators of $e_i$ with $\chi_{\ssc J,C}f_{\ol
rs}$. More precisely, we shall need the action of such commutators on
certain vectors $v_\l\in\VBL$. In this case we have:\par\ni {\sl Lemma
5.9.} Let $J=\{\ol p,\ldots,\ol1;1,\ldots,q\},\ C=(c_{\ol p},\ldots,c_{
\ol1};c_1,\ldots,c_q)$ with $p\le r\le m,\ q\le s\le n$. Let $v_\l\in
\VBL$ with weight $\l$ satisfying\vs{2pt}\nl\hs{50pt} $\matrix{c_{\ol r}
=\l(h_{\ol r})\mbox{ if }\ol p=\ol r,\hfill c_1-c_{\ol1}=\l(h_0),\hfill
c_s=\l(h_s)\mbox{ if }q=s,\vs{2pt}\cr c_{\ol i}\!-\!c_{\ol{i+1}}\!-\!1=
\l(h_{\ol i}){\rm\ if}\ \ol p<\ol i<0,\ \ \ c_i\!-\!c_{i+1}\!-\!1=\l(h_
i){\rm\ if}\ 0<i<q,\cr}$\hfill(5.10)\vs{2pt}\nl and\nl\hs{120pt} $f_{
\ol r,\ol{p+1}}v_\l=0\mbox{ \ if \ }p<r,\mbox{ \ \ \ \ \ }f_{q+1,s}v_\l=
0\mbox{ \ if \ }q<s.$\hfill(5.11) \nl Then for all $i\in I$,\vs{1pt}\nl
\hs{120pt} $[e_i,\chi_{\ssc J,C}f_{\ol rs}]v_\l=\cases{\chi_{\ssc J,C}
f_{\ol{r-1},s}v_\l&if \ $i=\ol r<\ol p$,\cr \chi_{\ssc J,C}f_{\ol r,s-1}
v_\l&if \ $i=s>q$,\cr \ \ 0&otherwise. \cr}$\hfill(5.12)\vs{3pt}\nl {\sl
Proof.} The first thing to note is that the only non-vanishing
commutators of $e_i$ with negative root vectors appearing in (5.8) are
the following:\vs{-6pt} $$\matrix{[e_{\ol i},f_{\ol a\ol i}]=f_{\ol a,
\ol{i+1}},\hfill&[e_{\ol i},f_{\ol i\ol i}]=h_{\ol i},\hfill&[e_{\ol i},
f_{\ol i\ol b}]=-f_{\ol{i-1},\ol b},\hfill&\mbox{ \ for\ }\ol a<i<b;\cr[
e_0,f_{\ol a0}]=f_{\ol a\ol1},\hfill&[e_0,f_{00}]=h_{0},\hfill&[e_0,f_
{0b}]=f_{1b},\hfill&\mbox{ \ for\ }\ol a<0<b;\cr[e_i,f_{ai}]=f_{a,i-1},
\hfill&[e_i,f_{ii}]=h_{i},\hfill&[e_i,f_{ib}]=-f_{i+1,b},\hfill&\mbox{
\ for\ }a<i<b,\cr}\eqno\matrix{(5.13{\rm a})\cr(5.13{\rm b})\cr(5.13{
\rm c})\cr}$$ \vs{-10pt}\nl Consider first $0<i<q$. The only
non-vanishing contributions to (5.12) arise from terms in (5.8) that do
not contain $c_ic_{i-1}$. These can be grouped together in sets of
three so that for any fixed $a<i<b$ the sum of each such set contains
the common factor $c_if_{ai}f_{i+1,b}-f_{a,i-1}f_{ii}f_{i+1,b}+c_{i+1}f_
{a,i-1}f_{ib}.$ Taking the commutator with $e_i$ and using (5.13c)
gives $c_if_{a,i-1}f_{i+1,b}-f_{a,i-1}h_if_{i+1,b}$ $-c_{i+1}f_{a,i-1}f_
{i+1,b}=f_{a,i-1}f_{i+1,b}(c_i-h_i-1-c_{i+1}),$ which acts to the right
on a sequence of $f_{i_x,i_{x+1}-1}$'s and $v_\l$. However, $[h_i,f_{i_
x,i_{x+1}-1}]=0$ since $i<b<i_x<i_{x+1}$ and $h_i$ acts finally on $v_
\l$ to give $\l(h_i)$. Thus all terms contain the common factor $c_i-\l(
h_i)-1-c_{i+1},$ which vanishes by virtue of our hypothesis (5.10). The
result for $\ol p<i<0$ is obtained in the same way. Similarly, $e_0$
commutes with all terms in (5.8) containing the product $c_{\ol1}c_1$.
The non-vanishing contributions to (5.12) can again be grouped into
sets of 3 terms such that for any fixed $\ol a${}$<$0$<${}$b$ the sum
of each such set contains the common factor $c_{\ol1}f_{\ol a0}f_{1b}$+$
f_{00}f_{\ol a\ol1}f_{1b}${}$-${}$c_1 f_{\ol a\ol1}f_{0b}.$ Taking the
commutator with $e_0$ and using (5.13b) gives $c_{\ol1}f_{\ol a\ol1}f_
{1,b}$+$h_0f_{\ol a\ol1}f_{1,b}${}$-${}$c_1f_{\ol a\ol1}f_{1b}$=$f_{\ol
a\ol1}f_{1,b}(c_{\ol1}$+$h_0${}$-${}$c_1)$. Moreover $h_0$ commutes
with everything else to its right to finally act on $v_\l$ giving $\l(h_
0)$. Thus all terms contain the common factor $c_{\ol1}$+$\l(h_0)${}$-
${}$c_1,$ again vanishes by (5.10).\par If $i$=$q$=$s$, $e_s$ commutes
with all terms in (5.8) other than those which can be paired so as to
give the common factor $-c_sf_{as}+f_{a,s-1}f_{ss}$ acting directly on $
v_\l$. Commutation with $e_s$ gives $-c_sf_{a,s-1}+f_{a,s-1}h_{ss}$
leading to the common factor $-c_s+\l(h_s),$ which vanishes. The result
for $i$=$\ol p$=$\ol r$ follows in the same way. If $i$=$q${}$<${}$s$, $
e_q$ commutes with every term in (5.8) other than those for which $i_
\ell$=$q$, but then $[e_q,f_{i_\ell s}]$=$[e_q,f_{qs}]$=$f_{q+1,s}$.
Thus every non-vanishing term contains the factor $f_{q+1,s}v_\l$ which
vanishes by (5.11). The story is the same for $i$=$\ol p${}$>${}$\ol r
$.\par For $\ol m${}$\le${}$i${}$<${}$\ol r$ or $\ol r${}$<${}$i${}$<
${}$\ol p$ or $q${}$<${}$i${}$<${}$s$ or $s${}$<${}$i${}$\le${}$n$ all
commutators with $e_i$ vanish since $i$ appears nowhere as a subscript
on any $f_{ab}$ appearing in (5.8). This leaves as non-vanishing only 2
special cases $i$=$s${}$>${}$q$ and $i$=$\ol r${}$<${}$\ol p$. In the
1st of these the only non-vanishing commutator of $e_s$ with terms in
(5.8) is $[e_s,f_{i_\ell s}]$=$f_{i_\ell,s-1}.$ This gives the 2nd case
of (5.12). Similarly the ${}$only non-vanishing commutator of $e_{\ol r}
$ with the terms in $(5.8)$ is $[e_{\ol r},f_{\ol r\ol{i_k}}]$=$f_{\ol
{r-1},\ol{i_k}}$, giving the 1st case of (5.12).\qed\par\ni Finally we
give a rather technical lemma which plays a crucial role in proving
results in \S6.\par\ni{\sl Lemma 5.10.} Given $\ol r,s,t$ with $\ol m
\le\ol r\le\ol1,\ 1\le s\le t\le n$. Let $J\subseteq\{\ol r,\cdots,\ol1,
1,\cdots,s\}$, $C=(c_{j_1},c_{j_2},...)\in\C^{\otimes\#J}$, $C(1)=(c_{j_
1}$+1,$c_{j_2}$+1,$...)$. Then \def \VSs{\vs{4pt}}\def\VSt{\vs
{5pt}}\VSs\nl\hs{100pt} $X_{\ssc J}\equiv\chi_{\ssc J,C}f_{\ol rs}\cdot
\chi_{\ssc J,C(1)}f_{\ol rt}+\chi_{\ssc J,C}f_{\ol rt}\cdot\chi_{\ssc J,
C(1)}f_{\ol rs}=0,$ \hfill(5.14)\VSt\nl where the 2nd term is obtained
from the 1st by interchanging $s$ and $t$. Taking $t=s$, we have \VSs\nl
\hs{130pt} $\chi_{\ssc J,C}f_{\ol rs}\cdot\chi_{\ssc J,C(1)}f_{\ol rs}=
0$\hfill(5.15)\VSs\nl {\sl Proof.} We prove this by induction on ${\sc
\#}J$. If ${\sc\#}J=0$, we immediately have $X_J=0$ since $f_{\ol rs}$,
$f_{\ol rt}$ anti-commute. Suppose now (5.14) holds for $J'$ with ${\sc
\#}J'<{\sc\#}J$. For $J$, suppose $J\cap I_2\ne\phi$ (the proof is
similar if $J\cap I_1\ne\phi$). Choose $j\in J$ to be the largest and
let $J'=J\bs\{j\}$. Let $C'$ and $C'(1)$ be respectively the element $C
$ and $C(1)$ corresponding to $J'$. Using (5.7) we have \vs{-4pt} $$
\mbox{the 1st summand of }X_J=\chi_{\ssc J',C'}(c_jf_{\ol rs}\!-\!f_{
\ol r,j-1}f_{js})\cdot\chi_{\ssc J',C'(1)}((c_j\!+\!1)f_{\ol rt}\!-\!f_
{\ol r,j-1}f_{jt}).\eqno(5.16)$$\vs{-14pt}\nl Now one may check the
validity of the following identities for all $\ol r\le i<j\le s$: \vs
{-6pt} $$\matrix{\chi_{i,c_i+1}(f_{\ol r,j-1}f_{js})=\chi_{i,c_i+1}f_{
\ol r,j-1}\cdot f_{js},\hfill\vs{3pt}\cr f_{js}\chi_{i,c_i+1}f_{\ol rt}=
\chi_{i,c_i+1}(f_{js}f_{\ol rt})=\chi_{i,c_i+1}f_{\ol rt}\cdot f_{js},
\hfill\vs{3pt}\cr f_{js}\chi_{i,c_i+1}(f_{\ol r,j-1}f_{jt})=\chi_{i,c_i+
1}f_{\ol rs}\cdot f_{jt}+\chi_{i,c_i+1}f_{\ol r,j-1}\cdot f_{js}f_{jt}.
\hfill\vs{-4pt}\cr}$$ Using these, (5.16) becomes\vs {-6pt} $$\matrix{c_
j(c_j+1)\chi_{\ssc J',C'}f_{\ol rs}\cdot\chi_{\ssc J',C'(1)}f_{\ol rt}-
c_j\chi_{\ssc J',C'}f_{\ol rs}\cdot\chi_{\ssc J',C'(1)}f_{\ol r,j-1}
\cdot f_{jt}\hfill\vs{2pt}\cr-(c_j+1)\chi_{\ssc J',C'}f_{\ol r,j-1}
\cdot\chi_{\ssc J',C'(1)}f_{\ol rt}\cdot f_{js}+\chi_{\ssc J',C'}f_{\ol
r,j-1}\cdot\chi_{\ssc J',C'(1)}f_{\ol rs}\cdot f_{jt}\hfill\vs{2pt}\cr+
\chi_{\ssc J',C'}f_{\ol r,j-1}\cdot\chi_{\ssc J',C'(1)}f_{\ol r,j-1}
\cdot f_{js}f_{jt}\hfill\vs{-5pt}\cr}$$ Denote these terms by $w_1,{\sc
\cdots},w_5$, and denote the corresponding terms for the 2nd summand of
$X_J$ by $w_6,{\sc\cdots},w_{10}$. Then $X_J=\sum^{10}_{k=1}w_k.$ By
the inductive hypothesis, we have $w_1+w_6=w_2+w_4+w_8=w_3+w_7+w_9=w_5=
w_{10}=0$. Hence $X_J=0$.\qed\par\ni The importance of this lemma lies
in the consequences which flow from the special case (5.15).\es {\bf
VI. PRIMITIVE VECTORS OF THE KAC-MODULE $\VBL$}\par\ni Let\,$\L$\,be a
dominant $r$-fold atypical weight of\,$G$\,with atypical roots\,\{$\g_1
$,${\sc\ldots}$,$\g_r$\}.\,In this section we first prove that to every
indecomposable unlinked code $\Si^c$ for $\L$ there corresponds a
primitive vector $v_{_\Si}$ of $\VBL$ having weight $\Si$. Then we
generalize the result to arbitrary unlinked codes. \vs {-8pt}\par As a
precursor to the proof we first restrict attention to those $\L$ for
which $\g_r=\a_{\ol mn}=\b_{1,n+1}$ and those codes $\Si^c$ for which $
D_\Si=SW_\L(r)$. It follows that the topmost and rightmost position $TR_
D$ of $D_\Si$ is $(1,n+1)$. Thus from (3.2),\vs {-7pt} $$A(\L)_{1,n+1}
\!=\!\rb{-5pt}{\mbox{$^{_{\ \sum}^{\sc\ \ 0}}_{k=\ol m}$}}a_k\!-\!\rb{-
5pt}{\mbox{$^{_{\ \sum}^{\sc\ \ n}}_{k=1}$}}a_k\!+\!m\!-\!n\!=\!0.\eqno(
6.1)$$\vs{-14pt}\nl It is convenient to introduce special labels for
some particular roots in $\wh D_\Si$. Let $x={\sc\#}\{j|(1,j)\in D_\Si
\}$, $y={\sc\#}\{i|i,n+1)\in D_\Si\}$ be respectively the number of
elements in the topmost row and rightmost column of $D_\Si$. Denote the
roots associated with the positions in the topmost row of $D_\Si$,
taken from right to left, by $\eta_1,{\sc\cdots},\eta_x$ and the roots
in the rightmost column of $D_\Si$, from top to bottom, by $\eta_1',{
\sc\cdots},\eta_y'$. Thus,\vs{-2pt} $\eta_1=\a_{\ol m,n},\cdots,\eta_x=
\a_{\ol m,n-x+1},\,\eta_1'=\a_{\ol m,n},\cdots,\eta_y'=\a_{\ol{m-y+1},n}
,$ with $\eta_1=\eta_1'=\b_{1,n+1}$. It should be noted the definitions
of $x$ and $y$ are such that:\vs{-2pt} $$a_{n-x+1}\!>\!0,a_{n-i+1}\!=
\!0,\,i=1,\cdots,x\!-\!1\mbox{ \ and \ }a_{\ol{m-y+1}}\!>\!0,a_{\ol{m-i+
1}}\!=\!0,\,i=1,\cdots,y\!-\!1.\eqno(6.2)$$\vs{-12pt}\nl All this may
be illustrated as follows in our $sl(6/5)$ Example 4.18 with $\L=[00020;
0;0210]$ encountered in \S\S3-4. $A(\L)$ was displayed in (3.3) with
codes enumerated in (3.4)\vs{-3pt} and chains in (4.3). $\Si^c$ is the
indecomposable unlinked code \rb{-8pt}{$^{^{^{\sc1\ 2}_{\sc2\ 5}}_{\sc5}
\,^{^{\sc3\ 4\ 5}_{\sc4\ 5}}_{\sc5}}$}. \vs{-5pt}Below we have set $D_
\Si$ alongside $A(\L)$, specifying all positions on the chain $SW_\L(
\ell)$ by an entry $\ell$ for $\ell=1,\cdots,5$. Here $x=2,y=4$ and in
the final array we identify the positions in $D_\Si$ associated with
roots $\eta_i$, $\eta'_i$:\nl\hs {100pt} $A(\L)=\left(^{^{^{\sc7\ 6\ 3\
1\ 0}_{\sc6\ 5\ 2\ 0\ \ol1}}_{\sc5\ 4\ 1\ \ol1\ \ol2}}_{^{^{\sc4\ 3\ 0\
\ol2\ \ol3}_{\sc1\ 0\ \ol3\ \ol5\ \ol6}}_{\sc0\ \ol1\ \ol4\ \ol6\ \ol7}}
\right),\ \ \ \ D_\Si=\left(^{^{^{\sc\cdot\,\ \cdot\,\ \cdot\,\ 5\ 5}_{
\sc\cdot\,\ \cdot\,\ 5\ 4\ 5_{\,}}}_{\sc\cdot\,\ \cdot\,\ 4\ 4\ 5_{\,}}}
_{^{\sc5\ 5\ 3\ 4\ 5_{\,}}_{^{\sc2\ 2\ 5\ \cdot\,\ \cdot_{\,}}_{\sc1\ 2
\ 5\ \cdot\,\ \cdot}}}\right)=\left(^{^{^{\sc\cdot\,\ \cdot\,\ \cdot\,\
\eta_2\ \eta_1}_{\sc\cdot\,\ \cdot\,\ *\ *\ \ \eta_2'}}_{\sc\cdot\,\
\cdot\,\ *\ *\ \ \eta_3'}}_{^{^{\sc*\ *\ *\ *\ \ \eta_4'}_{\sc*\ *\ *\
\cdot\,\ \ \cdot}}_{\sc*\ *\ *\ \cdot\,\ \ \cdot}}\right).$\hfill
(6.3)\nl In this example $\eta_1=\a_{\ol54},\eta_2=\a_{\ol53}$ and $
\eta_1'=\a_{\ol54},\eta_2'=\a_{\ol44},\eta_3'=\a_{\ol34},\eta_4'=\a_{
\ol24}$. Keeping this example in mind will help understand the proof
below.\par We shall always suppose 1$\le${}$x${}$\le${}$y$ as the
arguments for 1$\le${}$y${}$\le${}$x$ are entirely analogous. Set $$
\matrix{\L_0=\L,\hfill&\L_k=\L-\sum^k_{j=1}\eta_j,k=1,...,x,\hfill\cr v_
0(\L)=v_\L,\hfill&v_k(\L)=\chi_J^{(m/n-k+1)}f(\eta_k)v_{k-1}(\L),k=1,...
,x,\hfill\cr}\eqno(6.4)$$ with $J=\{\ol{m-y+1},\ldots,\ol1;1,\ldots,n-x+
1\}.$ Recall from Convention 2.2 that $\chi_J^{(m/n-k+1)}$ is the
operator $\chi_J$ defined for $sl(m+1,n-k+2)$ rather than for $G=sl(m+1/
n+1)$. \par\ni{\sl Lemma 6.1.} Let $1\le k\le x\le y$, then with
notation (5.6), \vs{2pt}\nl\hs{7ex}$\matrix{v_k(\L)=d_k(\L)v_{k-1}(\L),
\ d_k(\L)=\chi_{\ssc J,C}f_{\ol m,n-k+1},\hfill\cr c_j\!=\!c_j(\L)\!-
\!k\!+\!1,j\!\in\!J,\ c_j(\L)\!=\!\biggl\{^{\dis\SUM^j_{_{\sc\ell=\ol
{m-y+1}}}a_\ell\!+\!m\!+\!j,}_{\dis\SUM^{n-x+1}_{\ell=j}a_\ell\!+\!n\!-
\!j,}\,\ ^{^{\dis j\in\{\ol{m-y+1},\ldots,\ol1\},}}_{\dis j\in\{1,
\ldots,n-x+1\}.}\hfill\cr}${\hfill}\rb{2pt}{${^{^{^{^{\dis(6.5{\rm a})}}
}}_{_{_{\dis(6.5{\rm b})}}}}$} \vs{3pt}\nl {\sl Proof.} Since $v_{k-1}(
\L)$ has weight $\L_{k-1}$, it follows from definition (6.4) and
Corollary 5.7 that\vs{-7pt} $$v_k(\L)=\chi_{\ssc J,C}^{(m/n-k+1)}f_{\ol
m,n-k+1}v_{k-1}(\L),\ \ \ \ c_j=c_j^{(m/n-k+1)}(\L_{k-1}),\ j\in J.
\eqno(6.6)$$\vs{-17pt}\nl However, as can be seen from Corollary 5.8
with $r$=$m,s$=$n${}$-${}$k$+1, we have $\chi_{\ssc J,C}^{(m/n-k+1)}{
\sc\!}f_{\ol m,n-k+1}$ =$\chi_{\ssc J,C}{\sc\!}f_{\ol m,n-k+1}$, giving
as required (6.5a). Furthermore, note that $\L(h_i)$=$a_i,i${}$\in${}$I
$, the use of (6.2) in (5.9) immediately gives the 2nd equation of
(6.5b), so it remains to prove the 1st of (6.5b). However, it follows
from definitions (6.4) and (5.9) with ($m$/$n$) replaced by ($m$/$n${}$-
${}$k$+1) that \vs{-7pt} $$c_j^{(m/n-k+1)}(\L_{k-1})=\left\{\matrix{
\SUM^j_{_{\sc\ell=\ol{m-y+1}}}a_\ell+m-k+1+j,&j\in\{\ol{m-y+1},\ldots,
\ol1\},\vs{2pt}\cr\SUM^{n-x+1}_{\ell=j}a_\ell+n-k+1-j,\hfill&j\in\{1,
\ldots,n-x+1\}.\hfill\vs{-2pt}\cr}\right.\eqno(6.7)$$\vs {-12pt}\nl An
inspection of the 2nd of (6.5b) and (6.7) reveals that $c_j^{(m/n-k+1)}(
\L_{k-1})=c_j(\L)-k+1,\,j\in J.$ When used in (6.6) this completes the
proof of the 1st of (6.5b).\qed\par\ni {\sl Corollary 6.2.} Let $1\le k
\le x$, then $v_k(\L)\ne0$.\par\ni {\sl Proof.} Since $a_{\ol{m-y+1}}>0,
a_{n-x+1}>0$, (6.5b) implies $c_j(\L)>m+j\ge y-1\ge x-1,$ $j\in\{\ol{m-
y+1},\ldots,\ol1\}$; $c_j(\L)>n-j\ge x-1,$ $j\in\{1,\ldots,n-x+1\}$. It
follows that $c_j(\L)-k+1>0,\,j\in J,\,1\le k\le x$. However, $v_k(\L)=
{\sc\prod}_{\sc1\le i\le k,\,j\in J}(c_j(\L)-i+1)f(\eta_k)\cdots f(\eta_
1)v_{\ssc\L}+\cdots$, where the leading term is written in proper order
with the ordering (5.1). Thus $v_k(\L)\ne0$.\qed \par\ni Returning to $
sl(6/5)$ Example 4.18 with $\L=[0002;0;0210],\,x=2,y=4,\,J=\{\ol1;1,2,3
\}$. $\L_1=[\ol1002;0;021\ol1],\,\L_2=\L_x=[\ol2002;0;0200]$. The
coefficients $c_j(\L),\,j\in J$ are $(5;6,5,2)$ and $c_j$ associated
with $v_1(\L),\,v_2(\L)=v_x(\L)$ are $(5;6,5,2)$, $(4;5,4,1)$,
respectively. These are all non-zero, in accordance with Corollary 6.2.
A simpler example is Example 4.5, $\L=[1211;\ol1;10002],\,x=y=1,\,J=\{
\ol4,\ol3,\ol2,\ol1;1,2,3,4,5\}$, $\L_1=\L_x=[0211;\ol1;10001]$ and $c_
j(\L)$ associated with $v_1(\L)=v_x(\L)$ are $(1,4,6,8;7,5,4,3,2)$,
again all non-zero. To discuss the $G$-primitivity of $v_x(\L)$, it
should be noted that in the 1st of these 2 examples the weight $\L_x$
is not $G$-dominant, although the restriction of this weight to $G^{(m-
1/n)}=sl(5/5)$ is $G^{(m-1/n)}$-dominant. In contrast to this, in the
2nd example $\L_x$ is $G$-dominant. Guided by this distinction between
our 2 examples, it is convenient to deal first with a special case: \par
\ni{\sl Lemma 6.3.} $v_1(\L)$ is a $G$-primitive vector if $x=y=1$.\par
\ni {\sl Proof.} In this case, $J=\{\ol m,\ldots,\ol1;1,\ldots,n\}$ and
from (6.5) we have \nl\hs{12ex}$d_1(\L)\!=\!\chi_{\ssc J,C}f_{\ol mn},\
\ c_j\!=\!c_j(\L)\!=\!\cases{ \SUM^j_{_{\sc k=\ol m}}a_k +m + j,\!\!\!&
$j \in \{\ol m,\ldots,\ol1\}$, \cr \SUM^n_{k=j}a_k +n -j,\!\! \!&$j \in
\{1,\ldots,n\}$.\cr}$\nl It then follows that: \nl\cl{ $\matrix{\matrix
{c_{\ol m}=a_{\ol m},\hfill&c_n=a_n,\hfill\cr c_i\!-\!c_{i-1}\!-\!1=a_i
\hbox{ \ if \ }\ol m\!<\!i\!<\!0,\hfill&c_i\!-\!c_{i+1}\!-\!1=a_i\hbox{
\ if \ }0\!<\!i\!<\!n,\hfill\cr}\hfill\cr c_1-c_{\ol1}=\SUM^n_{k=1}a_k+
n-1-\SUM^{\ol1}_{k=\ol m}a_k-m+1=a_0,\hs{14pt}\hfill}$}\nl where, the
recovery of $a_0$ in the last case is a consequence of (6.1). It is
only here that the atypicality condition makes itself felt. Since $\L(h_
i)=a_i,\,i\in I$, it follows that $v_{\ssc\L}\in\VBL$ satisfies all
hypotheses (5.10) of Lemma 5.9 for $p=r=m,\,q=s=n$. The hypotheses
(5.11) are redundant, as are the first 2 cases of (5.12). Therefore, we
conclude from (5.12) that $[e_i,d_1(\L)]v_{\ssc\L}=0,\,i\in I$. Since $
v_{\ssc\L}$ is itself $G$-primitive, so $d_1(\L)v_{\ssc\L}$ is also $G
$-primitive.\qed\par\ni Prior to tackling other cases it is convenient
to introduce one further preliminary result:\par\ni{\sl Lemma 6.4.} $f_
{\ol m,\ol{m-y+2}}v_{\ssc\L}=0$ if $1<y$ and $f_{n-x+2,n-k+1}v_{\ssc\L}=
0$ if $1\le k<x.$ \par\ni{\sl Proof.} Let\,$I_{\ol y,x}$=$I_{\ol y}${}$
\cup${}$I_x,I_{\ol y}$=$\{{\sc\!}\ol m$,${\sc\ldots}$,$\ol{m\!{\sc\!}-{
\sc\!}\!y{\sc\!}\!+{\sc\!}\!2}{\sc\!}\},I_x$=${\sc\!}$\{${\sc\!}n${}$-
${}$x$+2,${\sc\ldots}$,$n{\sc\!}$\}.\,${\sc\!}$ Since\,$e_iv_{\ssc\L}
$=0\,and\,(6.2)\,gives $h_iv_{\ssc\L}$ =$0,i${}$\in${}$I_{\ol y,x}$,
consideration of algebra\, Span${\sc\!}$\{${\sc\!}e_i$,$f_i$,$h_i{\sc\!}
$\}\, implies\,that\,$f_iv_{\ssc\L}$=0. By (2.4),\, $f_{ij},i$,$j${}$
\in${}$I_{\ol y}$, $f_{k\ell},k$,$\ell${}$\in${}$I_x$\,
can\,respectively\,be\,expressed\,in\,terms\,of\,$f_i,i${}$\in${}$I_{
\ol y},f_j,j${}$\in${}$I_x$. The\,result\,then\,follows.\qed \par\ni
{\sl Lemma 6.5.} If $1<x\le y$, then $v_x(\L)$ is a $G^{(m-1/n)}=sl(m/n+
1)$ primitive vector.\vs{4pt}\nl\ni{\sl Proof.} In this case, by
(6.5b), we have \nl\hs{60pt} $\matrix{\matrix{c_{\ol{m-y+1}}=a_{\ol{m-y+
1}}+y-k,\hfill&c_{n-x+1}=a_{n-x+1}+x-k,\hfill\cr c_i\!-\!c_{i-1}\!-\!1=
a_i,\ \ol{m\!-\!y\!+\!1}\!<\!i\!<\!0,\hfill&c_{i}\!-\!c_{i+1}\!-\!1=a_i,
\ 0\!<\!i\!<\!n\!-\!x\!+\!1,\hfill\cr}\hfill\cr c_1-c_{\ol1}=\SUM^n_{k=
1}a_k+n-1-\SUM^{\ol1}_{k=\ol m}a_k-m+1=a_0,\hfill\cr}$\hfill (6.8)\vs
{2pt}\nl We are going to exploit Lemma 5.9 for $p$=$m${}$-${}$y$+1$,q
$=$n${}$-${}$x$+$1,r$=$m,s$=$n${}$-${}$k$+1 with 1$<${}$x${}$\le${}$y,$
1$\le${}$k${}$\le${}$x,v_\l$=$v_{k-1}(\L)$. It is necessary to check
that all hypotheses (5.10-11) are satisfied. Noted that $\l(h_i)$=$\L_
{k-1}(h_i)$=$a_i,\ol{m\!-\!1}${}$\le${}$i${}$\le${}$n${}$-${}$k$+1
implies $\l(h_i)$=$\L_{k-1}(h_i)$=$a_i$, $\ol{m\!-\!y\!+\!1}${}$\le${}$
i${}$\le${}$n${}$-${}$x$+1. It follows from (6.8) that the hypotheses
(5.10) are all satisfied unless either $i$=$\ol p$=$\ol{m\!-\!y\!+\!1}
$=$\ol r$= $\ol m$ with $y${}$\ne${}$k$, or $i$=$q$=$n${}$-${}$x$+1=$s
$=$n${}$-${}$k$+1 with $x${}$\ne${}$k$. Neither case can occur. Hence
(5.10) is satisfied. It remains to consider (5.11). From (2.4) $[f_{\ol
m,\ol{m-y+2}}$,$f_{ab}]$=0 unless $a$=$\ol{m\!-\!y\!+\!1}$ or $b$=$\ol
{m\!+\!1}$; $[f_{n-x+2,n-k+1}$,$f_{ab}]$=0 unless $a$=$n${}$-${}$k$+1
or $b$=$n${}$-${}$x$+1. However, the expansion of $d_i(\L)$ by means of
(5.8) involves only those $f_{ab}$ for which $a${}$\in$\{$\ol m$\}$\cup
$\{$\ol{m\!-\!y}$,${\sc\ldots}$,$n${}$-${}$x$+1\} and $b${}$\in$\{$\ol
{m\!-\!y\!+\!1}$,${\sc\ldots}$,$n${}$-${}$x$\}$\cup$\{$n${}$-${}$i$+1\}
 . This implies for $i$=1,${\sc\ldots}$,$k${}$-$1 we have $[f_{\ol m,\ol
{m-y+2}}$,$d_i(\L)]$=0,1$<${}$y$, $[f_{n-x+2,n-k+1}$,$d_i({\sc\!}\L{\sc
\!})]$=0,1$\le${}$k${}$<${}$x$. Since\,$v_{k-1}({\sc\!}\L{\sc\!})$=$d_
{k-1}({\sc\!}\L{\sc\!})${}${\sc\cdots}${}$d_1({\sc\!}\L{\sc\!})v_{\ssc
\L}$ it\,follows\,from\,Lemma\,6.4 that $f_{\ol m,\ol{m-y+2}}v_{k-1}(\L)
$=0, $f_{n-x+2,n-k+1}v_{k-1}(\L)$=0, confirming that (5.11) is
satisfied. Lemma 5.9 then implies that $[e_i,d_k(\L)]v_{k-1}(\L)=0$ for
$i\in I,1\le k\le x$ unless either $i=\ol m$ or $i=n-k+1,1\le k<x$. The
1st case does not concern us within $G^{(m-1/n)}$. The other cases
imply that since $v_x(\L)=g_xv_{\ssc\L},$ $g_x=d_x(\L)\cdots d_1(\L)$
we have $e_iv_x(\L)=g_xe_iv_{\ssc\L}=0$ unless $i=n-k+1,1\le k<x$. If $
i=n-k+1,1\le k<x$ we have\vs{3pt}\nl\cl{ $\matrix{e_{n-k+1}v_x(\L)\!\!
\!&=d_x(\L)\cdots d_{k+1}(\L)e_{n-k+1}d_k(\L)v_{k-1}(\L)\hfill\cr&=d_x(
\L)\cdots d_{k+1}(\L)[e_{n-k+1},d_k(\L)]v_{k-1}(\L)+d_x(\L)\cdots d_1(
\L)e_{n-k+1}v_{\ssc\L}\hfill\cr}$}\vs{3pt}\nl However, $e_{n-k+1}v_{
\ssc\L}=0$. Furthermore, the 2nd case of (5.12) and definitions (5.6\&
6.6) give $$d_{k+1}(\L)[e_{n-k+1},d_k(\L)]v_{k-1}(\L)=\pr{j\in J}\chi_
{j,c_j(\L)-k}f_{\ol m,n-k}\pr{j\in J}\chi_{j,c_j(\L)-k+1}f_{\ol m,n-k}v_
{k-1}(\L)=0,$$ where the final equality is a consequence of ($\ssc\!
$5.15$\ssc\!$) in Lemma$\ssc\!$ 5.10.$\ssc\!$ We conclude that $e_{n-k+
1}v_x({\ssc\!}\L{\ssc\!})$ $=0$ for $1\le k<x$, thereby completing the
proof that $v_x(\L)$ is $G^{(m-1/n)}$-primitive.\qed\par\ni {\bf
Theorem 6.6.} To any indecomposable unlinked code $\Si^c$ for $\L$,
there corresponds a primitive vector $v_{\ssc\Si}=g_{\ssc\Si}v_{\ssc\L}
$ of $\VBL$ with weight $\Si$ for some $g_{\ssc\Si}\in{\bf U}(G^-)$. \vs
{-3pt}\par\ni{\sl Proof.}\vs{-3pt} Suppose that the topmost and
rightmost position $TR_D$ of $D_\Si$ is $(m+1-m_\Si,n_\Si+1)$. \vs
{-1pt}By Theorem 4.17 there exists $t$ with $1\le t\le r$ such that $\g_
t=\a_{\ol{m_{_\Si}},n_{_\Si}},1\le m_{\ssc\Si}\le m,1\le n_{\ssc\Si}\le
n$. First we suppose\vs{-3pt} that $m_{\ssc\Si}=m,n_{\ssc\Si}=n$ so
that $t=r,\g_t=\a_{\ol m,n},TR_D=(1,n+1)$. We shall see later, \vs
{-3pt}by restriction from $G^{(m/n)}$ to $G^{(m_\Si/n_\Si)}$, we can
prove the theorem in general. Under this assumption, we are going to
prove it by induction on the depth $d={\sc\#}\wh D_\Si$. For $d=1$ for
which $D_\Si$ necessarily consists of the single position $(1,n+1)$ and
$x=y=1$. Thus, $\Si=\L-\a_{\ol m,n}$ is precisely the weight of the
vector $v_1(\L)=d_1(\L)v_{\ssc\L}$ which was shown to be $G$-primitive
in Lemma 6.3 and our Theorem 6.6 is satisfied by $v_{\ssc\Si}=v_1(\L)
$.\vs{-4pt}\par Let $d>1$. With notation (6.4), by Lemma 6.5, $v_x(\L)$
is $G^{(m-1/n)}$-primitive with weight $\L_x$. The atypicality matrix $
A(\L_x)$ is obtained from $A(\L)$ by subtracting $x$ from all entries
in the topmost row and adding 1 from each of the last $x$ columns. \vs
{-3pt}By removing the topmost row of $A(\L_x)$ we obtain the
atypicality matrix $A^{(m-1/n)}(\L_x)$ of $\L_x$ restricted to $G^{(m-1/
n)}$.\vs{-3pt} It is easy to check \vs {-3pt}that the $(2,n+1)$-th
entry in $A(\L_x)$ is always zero, so that $\b_{2,n+1}=\a_{\ol{m-1},n}$
is an atypical root for $\L_x$ restricted to $G^{(m-1/n)}$. \vs
{-2pt}Moreover, $\L_x$ is $r_x$-fold atypical with respect to $G^{(m-1/
n)}$ where\vs{-2pt} $r_x=r$ if $x<y$ or $r-1$ if $x=y$. Using results
in \S4, one can then see that there is an indecomposable unlinked code $
\Si_x^c$ for the restriction of $\L_x$ to $G^{(m-1/n)}$ for which\vs
{-9pt} $$\wh D_{\Si_x}=\wh D_\Si^{(m-1/n)}=\wh D_\Si\bs\{\eta_1,...,
\eta_x\}=\wh D(r_x)^{(m-1/n)},\eqno(6.9)$$\vs {-17pt}\nl
with\,the\,topmost\,and\,rightmost\,position\,of\,$D_\Si^{(m-1/n)}$
being\,$(1,n+1)$\,in\,$A(\L_x)^{(m-1/n)}$\,(but which is\,position (2,$
n$+1) in $A(\L_x)$). This may again be illustrated by $sl(6/5)$ example
with $\L$= $[00020;0;0210]$, $x$=$2,\L_x$=$[\ol10020;0;0200]$. On
restriction this yields $sl(5/5)$-dominant weight $\L_x$=$[0020;0;0200]
$ for which there exists an indecomposable unlinked code $\Si_x^c$
again\vs{-2pt} given by \rb{-10pt}{$^{^{^{\sc1\ 2}_{\sc2\ 5}}_{\sc5}\,^
{^{\sc3\ 4\ 5}_{\sc4\ 5}}_{\sc5}}$}. The corresponding atypicality
matrix $A(\L_x)$ and $D_{\Si_x}$ take the form ({\sl cf.}\,(6.3)):\vs
{-5pt}\nl\cl{ $A(\L_x)=\left(^{^{^{\sc5\ 4\ 1\ 0\ \ol1}_{\sc6\ 5\ 2\ 1\
0}}_{\sc5\ 4\ 1\ 0\ \ol1}}_{^{^{\sc4\ 3\ 0\ \ol1\ \ol2}_{\sc1\ 0\ \ol3\
\ol4\ \ol5}}_{\sc0\ \ol1\ \ol4\ \ol5\ \ol6}}\right),\ \ \ \ D_{\Si_x}=
\left(^{^{^{\sc\cdot\,\ \cdot\,\ \cdot\,\ \cdot\,\ \cdot}_{\sc\cdot\,\
\cdot\,\ 5\ 5\ 5_{\,}}}_{\sc\cdot\,\ \cdot\,\ 4\ 4\ 5_{\,}}}_{^{\sc5\ 5
\ 3\ 4\ 5_{\,}}_{^{\sc2\ 2\ 5\ \cdot\,\ \cdot_{\,}}_{\sc1\ 2\ 5\ \cdot
\,\ \cdot}}}\right)=\left(^{^{^{\sc\cdot\,\ \cdot\,\ \cdot\,\ \cdot\,\
\cdot}_{\sc\cdot\,\ \cdot\,\ *\ *\ \ \eta_2'}}_{\sc\cdot\,\ \cdot\,\ *\
*\ \ \eta_3'}}_{^{^{\sc*\ *\ *\ *\ \ \eta_4'}_{\sc*\ *\ *\ \cdot\,\ \
\cdot}}_{\sc*\ *\ *\ \cdot\,\ \ \cdot}}\right).$}\nl The positions of
the entries $*$ serve to specify the roots $\b${}$\in${}$\wh D_\Si^{(4/
4)}$.\,The result conforms precisely with (6.9) as can be seen by
comparison with the diagrams specifying the roots $\b\in D_\Si$ in
(6.3). Now let $U(\L_x)$ be the cyclic $G^{(m-1/n)}$-submodule ${\bf U}(
G^{(m-1/n)})v_x(\L)$, \vs{-2pt}which turns out to be isomorphic to \vs
{-3pt}the Kac-module $\ol V(\L_x)^{(m-1/n)}$ of $G^{(m-1/n)}$. Now the
depth of $\Si$ relative to $\L_x$ is given by $d_x={\sc\#}\wh D_\Si^{(m-
1/n)}={\sc\#}\wh D_\Si-x=d-x<d$, by induction hypothesis \vs {-3pt}we
see that there must exist\vs{-3pt} some $g_{\ssc\Si}^{(m-1/n)}\in{\bf U}
((G^-)^{(m-1/n)})$ such that $0\ne v_{_\Si}=g_{\ssc\Si}^{(m-1/n)}v_x(\L)
=g_{\ssc\Si}v_{\ssc\L}$ is $G^{(m-1/n)}$-primitive, with $g_{\ssc\Si}=g_
{\ssc\Si}^{(m-1/n)}g_x\in{\bf U}(G^-)$, $g_x=d_x(\L)\cdots d_1(\L)$.\vs
{-3pt}\par\ni {\sl Lemma 6.7.} $v_{_\Si}$ defined above is $G
$-primitive if $x<y$.\vs{-3pt}\par\ni{\sl Proof.} It remains to prove $
e_{\ol m}v_{\ssc\Si}=0.$ First,\vs{-3pt} in constructing $A(\Si)$ from $
A(\L)$ with $\Si=\L-\sum_{\b\in\wh D_\Si}\b,$ $A(\L)_{1,n+1}=0$, we
have $A(\Si)_{1,n+1}=y-x>0$.\vs{-3pt} Since $\Si$ is dominant we have $
A(\Si)_{1,i+1}>A(\Si)_{1,i+2}$ for $i=0,\ldots,n-1$. Hence, from (3.1a)
we have\vs{-9pt} $$A(\Si)_{1,i+1}=\la\Si+\rho|\a_{\ol mi}\ra>0\hbox{ \
\ for \ \ }i=0,\ldots,n.\eqno(6.10)$$\vs{-20pt}\nl \vs{-3pt}Second, we
can write $\Si=\L-\sum_{\b\in\wh D_\Si}\b=\L-\sum_{i=1}^d\b_j$ in the
way that $\b_i$ is an atypical root of $\Si_{i-1}=\L-\sum_{j=1}^{i-1}\b_
j$ (which is not necessarily $G$-dominant) in the sense that $\la\Si_{i-
1}+\rho|\b_i\ra=0$. Induction on $i$ gives $\la\Si_i+\rho|\Si_i+\rho\ra=
\la\L+\rho|\L+\rho\ra$. When $i=d$ we obtain $\la\Si+\rho|\Si+\rho\ra=
\la\L+\rho|\L+\rho\ra.$ Now suppose that $e_{\ol m}v_{\Si}\ne 0$. Let $
g$ be an element in ${\bf U}(G^+)$ with the largest possible weight $
\mu$ such that $u=gv_{\Si}\ne 0$. By this definition, $u$ is $G
$-primitive with weight $\Si'=\Si+\mu$. Using that $v_{\Si}$ is $G^{(m-
1/n)}$-primitive, $g$ can be chosen to be a sum of the form\vs{-6pt} $$
g=e_{\ol m,i_1}\cdots e_{\ol m,i_k},\mbox{ \ \ with \ \ }\mu=\SUM^{i_k}_
{j=i_1}\a_{\ol mj},\eqno(6.11)$$\vs{-16pt}\nl where $i_1\le\cdots\le i_
k$ and $i_j<\cdots<i_k$ if $i_j\ge 0$. Since $u$ has weight $\Si+\mu$,
Lemma 5.3 gives $\la\Si$+$\mu$+$\rho|\Si$+$\mu$+$\rho\ra=\la\L$+$\rho|
\L$+$\rho\ra,$ combined with $\la\Si$+$\rho|\Si$+$\rho\ra=\la\L$+$\rho|
\L$+$\rho\ra,$ we obtain\vs{-7pt} $$2\la\Si+\rho|\mu\ra+\la\mu|\mu\ra=0.
\eqno(6.12)$$\vs{-17pt}\nl Since $\mu$ has the form (6.11), denote it
by $\mu_k$. Induction on $k$ gives $\la\mu_k|\mu_k\ra=\la\mu_{k-1}|\mu_
{k-1}\ra+\la2\mu_{k-1}+\a_{\ol m,i_k}|\a_{\ol m,i_k}\ra\ge 0$, where
for the last inner product, we have made use of $\la\a_{\ol mi}|\a_{\ol
mj}\ra\ge 0$ for $i,j\in I$ (this can easily be proved by (2.2)). Also
we may prove as follows that $\la\Si+\rho|\a_{\ol mi}\ra>0$ for all $i
\in I$: if $i<0$, then $\a_{\ol mi}$ is a positive even root, so $\la
\Si|\a_{\ol mi}\ra\ge 0,$ $\la\rho|\a_{\ol mi}\ra>0$; if $i\ge 0$, then
$\a_{\ol mi}$ is a positive odd root and $\la\Si+\rho|\a_{\ol mi}\ra>0$
by (6.10). This proves that the {\it l.h.s.} of (6.12) is $>0$; this
contradiction proves that $e_{\ol m}v_{\Si}$ must be zero. Hence the
lemma follows.\qed\vs{-2pt}\par\ni{\sl Lemma 6.8.} Let $x$=$y$. If $v_{_
\Si}$ is not $G$-primitive then $e_{\ol mn}v_{_\Si}$ is primitive and $
e_{\ol mi}v_{_\Si}\ne0,\,i\le n$.\vs{-2pt}\par\ni{\sl Proof.} The proof
is very similar to that of Lemma 6.7. If $v_{_\Si}$ is not primitive,
we want to prove that $g$ in (6.11) must be $e_{\ol mn}$. The only
difference is that now $A(\Si)_{1,n+1}$=$y${}$-${}$x$=0. Thus (6.10)
must be replaced by the statement $\la\Si$+$\rho|\a_{\ol mi}\ra${}$>${}$
0,i$=0,$\sc\ldots$,$n$-1,$\la\Si$+$\rho|\a_{\ol mn}\ra$=0. Hence, by
the same argument as before, for our hypothesized $G$-primitive vector $
u$=$gv_{_\Si},g$ must have weight $\mu$=$\a_{\ol mn}$, since this is
the only possible solution of (6.12). It follows that $g$=$e_{\ol mn}$.
By our choice of $g$, $e_{\ol mn}v_{\ssc\Si}\ne 0$ and so $e_{\ol mi}v_
{\ssc\Si}${}$\ne$0 since $e_{\ol mn}v_{\ssc\Si}$=$-${}$e_{in}e_{\ol mi}
v_{\ssc\Si}$ for $i<\!n$.\qed \vs{-2pt}\par\ni By Lemma 6.7, if $x<y$,
the proof of the theorem is then completed. So, let $x$=$y$. If\vs
{-2pt} $v_{_\Si}$ is $G$-primitive, the proof is also completed.
Suppose now $v_{_\Si}$ is not $G$-primitive. First note that \vs{-5pt}
$$\matrix{\wh D_{\ssc\Si}\!\!\!\!&=\!\!\!\!&\wh D_{\ssc\Si}^{(m-1/n)}{
\sc\cup}\{\eta_1,{\sc\ldots},\eta_x\}=\wh D_{\ssc\Si}^{(m/n-1)}{\sc\cup}
\{\eta_1',{\sc\ldots},\eta_x'\}\hfill\cr&=\!\!\!\!&\wh D_{\ssc\Si}^{(m-
1/n-1)}{\sc\cup}\{\a_{\ol mn}\}{\sc\cup}\{\eta_2,{\sc\ldots},\eta_x\}{
\sc\cup}\{\eta_2',{\sc\ldots},\eta_x'\}.\hfill\cr}\eqno(6.13)$$ \vs
{-12pt}\nl Then we see that $v_{_\Si}$ defined above can be written in
the form \vs{-6pt} $$v_{_\Si}=g_{_\Si}v_{\ssc\L}=g_{_\Si}^{(m-1/n)}v_x(
\L)=g_{_\Si}^{(m-1/n)}g_xv_{\ssc\L}=g_{_\Si}^{(m-1/n-1)}g'_{x-1}g_xv_{
\ssc\L},\eqno(6.14)$$ \vs{-15pt}\nl \vs{-3pt}where, quite generally the
weight ${\rm wt}(g_{_\Si}^{(r/s)})$=$-${}$\sum_{\a\in\wh D_\Si^{(r/s)}}
\a$ and ${\rm wt}(g_x)$=$-${}$\sum_{i=1}^x\eta_i,$ ${\rm wt}(g'_{x-1})
$= $-${}$\sum_{i=2}^x\eta_i'$. Furthermore, by induction on ${\sc\#}D_
\Si$, (6.13-14) tells that we can decompose\vs {-4pt} $$\matrix{D_\Si\!=
\!{\sc\cup}_{i=0}^{i_{_\Si}}D^{(i)},\hfill&D^{(1)}\!=\!\{\eta'_2,{\sc...
},\eta'_x\},\hfill&D^{(0)}\!=\!\{\eta_1,{\sc...},\eta_x\},\hfill\cr g_{_
\Si}\!=\!g^{(i_{_\Si})}{\sc\cdots}g^{(1)}g^{(0)},\hfill&g_{_\Si}^{(m-1/
n-1)}\!=\!g^{(i_{_\Si})}{\sc\cdots}g^{(2)},\hfill&g^{(1)}\!=\!g'_{x-1},
\ g^{(0)}\!=\!g_x,\cr}\eqno^{\dis(6.15{\rm a})}_{\dis(6.15{\rm b})}$$\vs
{-11pt}\nl \vs{-4pt}for some $i_{_\Si}$ such\,that$\,{\rm wt}(g^{(i)})
$=$-${}${\sc\sum}_{\a\in\wh D^{(i)}}\a,\,$where $D^{(i)}\,
$consists\,all\,positions of either the topmost row or the rightmost \vs
{-4pt}column in $D_\Si${}$\bs${}${\sc\cup}_{j=0}^{i-1}D^{(j)}$. Now
suppose the\vs{-3pt} position $P_x$=($x$,$n$+1$-${}$x$), which is$\,
$clearly\,in$\,D_\Si,\,$ belongs\,to$\,D_\Si^{(j_{_\Si})}\,$for\,some$
\,j_{_\Si}$: $\!2${}$\le${}$j_{_\Si}${}$\le${}$i_{_\Si},\,
$then\,we\,can\,write\vs{-6pt} $$\matrix{g_{_\Si}\!=\!g_{_\Si}^{(2)}g_{_
\Si}^{(1)},\hfill&g_{_\Si}^{(2)}\!=\!g^{(i_{_\Si})}{\sc\cdots}g^{(j_{_
\Si}+1)},\hfill&g_{_\Si}^{(1)}\!=\!g^{(j_{_\Si})}{\sc\cdots}g^{(0)},
\hfill\vs{1pt}\cr D_\Si\!=\!D_\Si^{(2)}{\sc\cup}D_\Si^{(1)},\hfill&D_
\Si^{(2)}\!=\!{\sc\cup}^{i_{_\Si}}_{i=j_{_\Si}+1}D^{(i)},\hfill&D_\Si^{(
1)}\!=\!{\sc\cup}^{j_{_\Si}}_{i=0}D^{(i)},\hfill\cr}\eqno^{\dis(6.15{
\rm c})}_{\dis(6.15{\rm d})}$$\vs{-16pt}\nl with \ \vs{-7pt} ${\rm wt}(
g_{_\Si}^{(i)})=\sum_{\a\in D_\Si^{(i)}}\a,$ $i=1,2.$ \ As an example,
consider $sl(5/6)$ with $\L=[0011;1;00200]$, $\Si^c$=\rb{-3pt}{$\,^{^{
\sc 1\,2\,3\,4}_{\sc3\,3\,4}}_{\sc4\,4}$}.\,Below\,we have\vs{-5pt} set$
\,D_\Si\,$alongside$\,A(\L),\,$specifying all positions in$\,D^{(i)}\,
$by an entry$\,i\,$for $i$=0,${\sc\ldots}$,8.\,\vs{-2pt}Here$\,x$=$y$=$
3,i_{_\Si}$=$8,j_{_\Si}
$=4.\,In\,the\,final\,array\,we\,identify\,the\,positions\,in$\,D^{(1)}_
\Si\,$by\,1.\vs{-8pt} $$A(\L)\!=\!\left(\rb{-3pt}{\mbox{$^{^{^{ \sc
7\,6\,5\,2\,1\,0}_{ \sc6\,5\,4\,1\,0\,\ol1}}_{ \sc5\,4\,3\,0\,\ol1\,\ol
2}}_{^{ \sc3\,2\,1\,\ol2\,\ol3\,\ol4}_{ \sc1\,0\,\ol1\,\ol4\,\ol5\,\ol
6}} $}}\right),\ D_\Si\!=\!\left(\rb{-3pt}{\mbox{$^{^{^{ \sc.\ .\ .\,
0\,0\,0}_{ \sc.\ .\ .\, 2\,2\,1}}_{ \sc4\,4\,4\,4\,3\,1}}_{^{ \sc
7\,7\,6\,5\,3\,.}_{ \sc8\,8\,6\,5\,.\ .}} $}}\right)\!=\!\left(\rb{-3pt}
{\mbox{$^{^{^{ \sc.\ .\ .\, 1\,1\,1}_{ \sc.\ .\ .\, 1\,1\,1}}_{ \sc
1\,1\,1\,1\,1\,1}}_{^{ \sc*\,*\,*\,*\,1\,.}_{ \sc*\,*\,*\,*\,.\ .}} $}}
\right).\eqno(6.16)$$\vs{-16pt}\nl From this example, it it easy to
obtain the following result. \vs{-6pt}\par\ni{\sl Lemma 6.9. }\vs{-3pt}
$\!$(i) $\!$All\,positions\,above\,and\,to\,the\,right\,of\,$P_x
$\,are\,in\,$D_\Si^{(1)}$,\,i.e.,\,($i$,$j$)$\in${}$D_\Si^{(1)},1${}$
\le${}$i${}$\le${}$x,$ $n$+1 $-${}$x${}$\le${}$j${}$\le${}$n
$+1;\,none\,position\,below\,and\,to\,the\,left\,of\,$P_x$\,is\,in\,$D_
\Si^{(1)}$,\,i.e.,\,($i$,$j$)$\not\in${}$D_\Si^{(1)},i${}$>${}$x$,$j${}$
<${}$n$+1$-${}$x$; \vs{-3pt}\nl \hs{3ex}(ii)(a)\vs{-3pt} If ($i$,$j$)$
\in${}$D_\Si^{(1)},i${}$\le${}$x$ then ($i$,$j'$)$\in${}$D_\Si^{(1)}$
for all $j'$: $\!j${}$\le${}$j'${}$\le${}$n$+1. (b) If ($i$,$j$)$\in
${}$D_\Si^{(1)},n$+1$-${}$x${}$\le${}$j$ then ($i'$,$j$)$\in${}$D_\Si^{(
1)}$ for all $i'$: $\!1${}$\le${}$i'${}$\le${}$i$. \qed\vs{-3pt}\par \ni
\vs{-3pt}Now remember that $x$=$y$, we can construct another vector $
\wt v_x(\L)$ if we start from the last column of $D_\Si$ instead of
the\vs{-3pt} 1st row, such that $\wt v_x(\L)$=$\wt g_xv_{\ssc\L}$,
which has similar properties to those of $v_x(\L)$ but where ${\rm wt}(
\wt g_x)$=$-${}$\sum_{i=1}^x\eta'_i$. Then from $\wt v_x(\L)$, we can
construct a vector\vs{-8pt} $$\wt v_{_\Si}\!=\!\wt g_{_\Si}v_{\ssc\L}\!=
\!\wt g_{_\Si}^{(m/n-1)}\wt v_x(\L)\!=\!\wt g_{_\Si}^{(m/n-1)}\wt g_xv_
{\ssc\L}\!=\!\wt g_{_\Si}^{(m-1/n-1)}\wt g_{x-1}'\wt g_xv_{\ssc\L},
\eqno(6.17)$$\vs{-18pt}\nl which\,is$\,G^{(m/n-1)}$-primitive\,and ${
\rm wt}(\wt g_{x-1}')$=$-${}$\sum_{i=2}^x\eta_i$. Note\,that\,both $g^{(
m-1/n-1)}_{_\Si}$,\vs{-4pt} $\wt g^{(m-1/n-1)}_{_\Si}$ are
the\,element\,$g_{{\ssc\Si}^{(m-1/n-1)}_1}$\vs{-6pt}\,defined\,for$\,
$the\,weight $\Si_1$=$\L${}$-${}$\sum_{i=1}^x\eta_i${}$-${}$\sum_{i=2}^
x\eta'_i$ restricted\,to $G^{(m-1/n-1)}$, so,\,by\,induction\,on$\,m,n,
\,$we\,can\,suppose$\,g^{(m-1/n-1)}_{_\Si}$=$\wt g^{(m-1/n-1)}_{_\Si}\,
$and\,so\,by\,(6.13-15),\,we\,have\vs{-12pt} $$v_{\ssc\Si}\!=\!g^{(2)}_
{\ssc\Si}g^{(1)}_{\ssc\Si}v_{\ssc\L},\mbox{ \ and \ }\wt v_{\ssc\Si}\!=
\!g^{(2)}_{\ssc\Si}\wt g^{(1)}_{\ssc\Si}v_{\ssc\L},\,\wt g^{(1)}_{\ssc
\Si}\!=\!g^{(j_{_\Si})}{\sc\cdots}g^{(2)}\wt g_{x-1}'\wt g_x\eqno(6.18)
$$\vs{-18pt}\nl Now\vs{-3pt} if \rb{1pt}{$\wt v_{_\Si}$} is $G
$-primitive, then the proof is again completed, or by analogy with
Lemma 6.8, $e_{\ol mn}\wt v_{_\Si}$ is $G$-primitive. Thus, we can
suppose\vs{-11pt} $$e_{\ol m}v_{_\Si}\!\ne\!0\!\ne e_n\wt v_{_\Si},
\mbox{ \ \ but both $e_{\ol mn}v_{_\Si}$ and $e_{\ol mn}\wt v_{_\Si}$
are $G$-primitive.}\eqno(6.19)$$\vs{-14pt}\nl {\sl Lemma 6.10.} Let $v_
{_\Si}$ and $\wt v_{_\Si}$ be as in (6.19). Then $e_{\ol mn}v_{_\Si}=e_
{\ol mn}\wt v_{_\Si}$ (up to a non-zero scalar).\vs{-3pt} \par\ni{\sl
Proof.} By (6.18), we have\vs{-9pt}\par\ni \hs{100pt}$e_{\ol mn}v_{_\Si}
\!=\!g_{_\Si}^{(2)}v_2,\ v_2\!=\!u_2v_{_\L}\mbox{ \ \ and \ \ }e_{\ol
mn}\wt v_{_\Si}\!=\!g_{_\Si}^{(2)}\wt v_2,\ \wt v_2\!=\!\wt u_2v_{_\L},
$\hfill(6.20)\vs{-4pt}\nl \vs{-4pt}such\,that$\,u_2$=$[e_{\ol mn}$,$g_{_
\Si}^{(1)}],\wt u_2$=$[e_{\ol mn}$,$\wt g_{_\Si}^{(1)}]$
have\,the\,weight $-${}$\SUM_{\b\in\wh D_\Si^o}\b,D_\Si^o$=$D_\Si^{(1)}
\bs\{a_{\ol mn}\}.$ By Lemma 6.11\,below,\,we see that the only
possible prime term of $v_2$ is $b_1y_1v_{\ssc\L},y_1\!\in\!\C$.\vs
{-1pt} If $y_1$=0 then $v_2$ has no prime term,\vs{-1pt} and by Lemma
5.1(i), $e_{\ol mn}v_{\ssc\Si}$ has no prime term, which contradicts
with Lemma 5.2(i). Therefore $y_1${}$\ne$0. Similarly $\wt v_2$ has one
prime term $b_1\wt y_1v_{\ssc\L}$. By rescaling, we can suppose $y_1$=$
\wt y_1$ and then by Lemma 5.1(ii) and Lemma 5.2(ii), we obtain $e_{\ol
mn}v_{\ssc\Si}$=$e_{\ol mn}\wt v_{\ssc\Si}$.\qed \vs{-1pt}\par\ni{\sl
Lemma 6.11.} There is a unique $b_1={\sc\prod}_{\b\in\wh D_\Si^o}f(\b)
\in B$ with weight $-\sum_{\b\in\wh D_\Si^o}\b$.\vs{-1pt}\par \ni{\sl
Proof.}\hfill Using (6.16) as an example will help us to understand the
proof. Suppose there exists\par\ni another $b_2={\sc\prod}_{\b\in\wh D_
1}f(\b)\in B,\,D_1\subset D$ with the same weight $-\g$ so that $\g=
\sum_{\b\in\wh D_1}\b=\sum_{\b\in\wh D_\Si^o}\b$. Let $\g=\sum_{i=\ol m}
 ^na_i\a_i$ with the coefficients $a_i$, then we have \vs{-7pt} $$
\matrix{{\sc\#}\{i|(i,j)\in D_1\}={\sc\#}\{i|(i,j)\in D_\Si^o\}=a_{j-1}-
a_j,\hfill&1\le j\le n+1,\hfill\cr{\sc\#}\{j|(i,j)\in D_1\}={\sc\#}\{j|(
i,j)\in D_\Si^o\}=a_{\ol{m+1-i}}-a_{\ol{m+2-i}},\hfill&1\le i\le m+1,
\hfill\cr}\eqno\matrix{(6.21{\rm a})\cr(6.21{\rm b})\cr}$$ \vs{-12pt}\nl
where, if $j$=$n$+1 or $i$=1 we suppose $a_{n+1}$=$a_{\ol{m+1}}$=0. It
suffices to prove that the solution in (6.21) is $D_1$=$D_\Si^o$. If
not, suppose $(i_0$,$j_0$)$\in${}$D_1\bs${}$D_\Si^o$ with $j_0$ being
the smallest. First, suppose $(i_0$,$j_0$)=(1,$n$+1). Then (6.21a) in
the case $j$=$n$+1 tells us that there exists $(i_1$,$n$+1)$\in${}$D_
\Si^o\bs${}$D_1$ with $i_1${}$\ne${}$i_0$ and (6.21b) in the case $i$=$
i_1$ tells us that there exists $(i_1$,$j_1$)$\in${}$D_1\bs${}$D_\Si^o$
with $j_1${}$\ne${}$n$+1 and so, $j_1${}$<${}$n$+1=$j_0$, which
contradicts with the choice of $j_0$ being the smallest. Second,
suppose $(i_0$,$j_0$)$\ne$(1,$n$+1). Then by Lemma 6.9(i), we have $i_0
${}$>${}$x$ or $j_0${}$<${}$n${}$-${}$x$+1. Suppose $j_0${}$<${}$n${}$-
${}$x$+1 (the other case is similar). Then using (6.21a) in the case $j
$=$j_0$, there exists some $(i_1$,$j_0$)$\in${}$D_\Si^o\bs${}$D_1$. By
Lemma 6.9(i) we have $i_1${}$\le${}$x$. Also $i_1${}$\ne$1 since $D$,
and so, $D_\Si^o$, does not contain position (1,$j_0)$. Now consider
(6.21b) for $i$=$i_1$, from $(i_1$,$j_0$)$\in${}$D_\Si^o\bs${}$D_1$,
there exists some $j_1$ such that $(i_1$,$j_1$)$\in${}$D_1\bs${}$D_\Si^
o$; however, if $j_0${}$\le${}$j_1$, by Lemma 6.9(ii)(a) we would have $
(i_1$,$j_1$)$\in${}$D_\Si^o$, thus, $j_1${}$<${}$j_0$, again
contradicting with the choice of $j_0$ being the smallest.\vs{-3pt}\qed
\par\ni Now we can complete the proof of Theorem 6.6 as follows. By
Lemma 6.10, we can suppose $e_{\ol mn}(v_{_\Si}${}$-${}$\wt v_{_\Si})
$=0. Since $e_{\ol m}(v_{_\Si}${}$-${}$\wt v_{_\Si})$=$e_{\ol m}v_{_\Si}
${}$\ne$0, let $k$ be the largest such that $v_\l${}$\equiv${}$e_{\ol
mk}(v_{_\Si}${}$-${}$\wt v_{_\Si})${}$\ne$0 with $k${}$<${}$n$. As $\wt
v_{_\Si}$ is $G^{(m/n-1)}$-primitive, we have $v_\l$=$e_{\ol mk}v_{_\Si}
$ with weight $\l$=$\a_{\ol mk}$+$\Si.$ Applying $e_i$ to $v_\l$ and
(2.4) gives $e_iv_\l\!=\!0,\,i>\ol m$. We also have $e_{\ol m}v_\l\!=
\!0$: if not, again similar to the arguments after (6.11), let $g_1${}$
\in${}${\bf U}(G^+)$ with largest possible weight $\mu_1$ such that $u
$=$g_1v_\l$=$gv_{_\Si}${}$\ne$ 0, where now $g$=$g_1e_{\ol mk}$ with
weight $\mu\!=\!\mu_1$+$\a_{\ol mk}$, then as $\la\Si$+$\rho|\a_{\ol mk}
\ra\!>\!0$, we could find no solution for $\mu$ (or $\mu_1$) in (6.12).
Thus we have in fact proved that $v_\l$ is $G$-primitive, which is not
possible since $e_{kn}v_\l\!=\!e_{\ol mn}v_{\ssc\Si}\!\ne\!0$. The
contradiction shows that the assumption (6.19) is wrong, so,\,either$
\,v_{_\Si}\,$or$\,\wt v_{_\Si}\,$must\,be$\,G$-primitive,$\,
$proving\,Theorem\,6.6\,in\,the\,case$\,TR_D$=(1,$n$+1).\vs{-3pt}\par
For $TR_D\!=\!(m\!+\!1\!-\!m_{\ssc\Si},n_{\ssc\Si}+\!1)$, let $G'=G^{(m_
{_\Si}/n_{_\Si})}$ and let $U^{(m_{_\Si}/n_{_\Si})}$ be the $G'
$-submodule of $\VBL\,$generated by$\,v_{\ssc\L}\,$isomorphic to$\,\ol
V(\L)^{(m_{_\Si}/n_{_\Si})}.\,$Let$\,\Si^{(m_{_\Si}/n_{_\Si})}\,
$correspond to an indecomposable unlinked code $\Si^{c(m_{_\Si}/n_{_\Si}
)}$ of $\L$ restricted to $G'$. By construction the topmost and
rightmost position of $D_\Si^{(m_{_\Si}/n_{_\Si})}$ in $A(\L)^{(m_{_\Si}
/n_{_\Si})}$ is $(1,n_{\ssc\Si}+1)$. As just proved, there is a $G'
$-primitive vector $v_{_\Si}=g_{_\Si}v_{\ssc\L}$ corresponding to the
code $\Si^{c(m_{_\Si}/n_{_\Si})}$ with $g_{_\Si}\in{\bf U}(G'^-)$,
which commutes with $e_i,\,i\in\{\ol m,\ol{m-1},\cdots,\ol{m_{{\ssc\Si}}
+1},n_{\ssc\Si}+1,\cdots,n\}$. Hence $v_{_\Si}$ is also $G$-primitive
corresponding to the code $\Si^c$. This completes the proof of Theorem
6.6 in general.\qed\vs{-3pt}\par\ni{\bf Theorem 6.12.} To any unlinked
code $\Si^c$ for $\L$, there corresponds a primitive vector $v_{\ssc\Si}
=g_{\ssc\Si}v_{\ssc\L}$ of $\VBL$ with weight $\Si$ for some $g_{\ssc
\Si}\in{\bf U}(G^-)$. \vs{-3pt}\par\ni {\sl Proof.} Suppose $\Si^c=\Si_
1^c\cdots\Si_k^c$ with all $\Si_i^c$ indecomposable unlinked codes. The
proof is covered by Theorem 6.6 if $k=1$. Let $k>1$. By Definition
3.12, we see that $\Si_k^c$ is an indecomposable \,code \,of \,$\L,\,$
thus corresponding to \,a \,primitive \,vector $\,v_{\ssc\Si_k}=g_{\ssc
\Si_k}v_{\ssc\L}$.\hfill We \,also\par\ni see that $\Si_o^c=\Si_1^c
\cdots\Si_{k-1}^c$ is an unlinked code for the highest weight vector $v^
*_{\ssc\Si_k}$ of $\ol V(\Si_k)$ with weight $\Si_k$. By induction on $
k$, there exists a primitive vector $v^*_{\ssc\Si_o}=g_{\ssc\Si_o}v^*_{
\ssc\Si_k}$ in $\ol V(\Si_k)$. Let $v_{_\Si}=g_{\ssc\Si_o}g_{\ssc\Si_k}
v_{\ssc\L}\in\VBL$ be the image of $v^*_{\ssc\Si_o}$ under the
homomorphism $\ol V(\Si_k)\rar\VBL$: $v^*_{\ssc\Si_k}\rar v_{\ssc\Si_k}
$, $gv^*_{\ssc\Si_k}\rar gv_{\ssc\Si_k}$, $g\in{\bf U}(G)$. One can
check that $v_{_\Si}$ is nonzero as its leading term up to a non-zero
scalar is ${\sc\prod}_{\b\in\wh D_\Si}f(-\b)v_{\ssc\L}$, hence it is a
primitive vector corresponding to $\Si$.\qed \par\ \vs{-7pt}\par\ni {\bf
ACKNOWLEDGEMENTS}\vs{-2pt}\par\ni We would like to thank J. Van der
Jeugt for many useful discussions. YCS would like to thank Queen Mary
and Westfield College and Concordia University for hospitality during
the Academic year 1990-91 and 1991-93 when most part of this work was
conducted and in particular, he thinks Profs. C.J.Cummins and J.McKay
for various helps. \vs{-7pt}\par\ \par\ni{\bf REFERENCES}\vs{-2pt}\par
\def\r{\par\ni\hi2.5ex\ha1} \ni\hi2.5ex\ha1 $^{\ 1}$ J. Germoni,
``Indecomposable representations of special linear {\parskip .02 truein
\baselineskip 2pt \lineskip 2pt Lie superalgebras,'' {\sl J. Alg.} {\bf
209} (1998), 367--401. \r $^{\ 2}$ J.W.B. Hughes, R.C. King and J. Van
der Jeugt, ``On the composition factors of Kac modules for the Lie
superalgebra $sl(m/n)$,'' {\sl J. Math. Phys.}, {\bf 33}(1992),
470--491.\r $^{\ 3}$ V.G. Kac, ``Classification of simple Lie
superalgebras,'' {\sl Funct. Anal. Appl.}, {\bf9} (1975), 263--265.\r $^
{\ 4}$ V.G. Kac, ``Lie superalgebras,'' {\sl Adv. in Math.} {\bf
26}(1977), 8--96.\r $^{\ 5}$ V.G. Kac, {\sl Characters of typical
representations of classical Lie superalgebras,'' {\sl Comm. Alg.} {\bf
5} (1977), 889--897.\r $^{\ 6}$ V.G. Kac, ``Representations of
classical Lie superalgebras,'' in {\sl Lecture Notes in Math.}, {\bf
676} (1978), 597--626.\r $^{\ 7}$ V.G. Kac, {\sl Infinite dimensional
Lie algebras}, 3rd ed., Cambridge University Press, 1990.\r $^{\ 8}$
T.D.Palev, ``Finite-dimensional representations of the special linear
Lie superalgebra $sl(1,n)$ II. Nontypical representations,'' {\sl J.
Math. Phys.} {\bf29} (1988), 2589-2598\r $^{\ 9}$ T.D.Palev,
``Essentially typical representations of the Lie superalgebras $gl(n/m)
$ in a Gel'fand-Zetlin basis,'' {\sl Funct. Anal. Appl.} {\bf 23}
(1989), 141-142 (English translation).\r $^{10}$ M. Scheunert, {\sl The
theory of Lie superalgebras,} Springer, Berlin, 1979.\r $^{11}$ V.
Serganova, ``Kazhdan-Luszitig polynomials and character formula for the
Lie superalgebra $gl$,'' {\sl Selecta Mathematica} {\bf2} (1996),
607--651.\r $^{12}$ Y. Su, ``Classification of finite dimensional
modules of the Lie superalgebra $sl(2/1)$,'' {\sl Comm. Alg.} {\bf20}
(1992), 3259--3277.\r $^{13}$ Y. Su, ``Primitive vectors of Kac-modules
of the Lie superalgebras $sl(m/n)$: II,'' in preparation.\r $^{14}$ J.
Van der Jeugt, Characters and composition factor multiplicities for the
Lie superalgebra $gl(m/s)$,'' {\sl Lett. Math. Phys.} {\bf47} (1999),
49--61.\r $^{15}$ J. Van der Jeugt, J.W.B. Hughes, R.C. King and J.
Thierry-Mieg, ``Character formulas for irreducible modules of the Lie
superalgebra $sl(m/n)$,'' {\sl J. Math. Phys.}, {\bf 31}(1990),
2278--2304.\r $^{16}$ J. Van der Jeugt, J.W.B. Hughes, R.C. King and J.
Thierry-Mieg, ``A character formula for singly atypical modules of the
Lie superalgebra $sl(m/n)$,'' {\sl Comm. Alg.}, {\bf 19}(1991),
199--222. }\end{document}